\documentclass[11pt]{article}
\usepackage[english]{babel}
\usepackage[letterpaper,margin=2cm]{geometry}
\usepackage[T1]{fontenc}
\usepackage[utf8]{inputenc}
\usepackage{lmodern}
\usepackage{amsmath,amssymb,amsthm,bm}
\usepackage{graphicx,booktabs,float,caption,subcaption}
\usepackage[dvipsnames]{xcolor}
\usepackage[numbers]{natbib}
\usepackage[colorlinks=true,allcolors=MidnightBlue]{hyperref}
\hypersetup{pdftitle={A High-Order Geometry-Conforming Selective Discontinuous Galerkin Method and Its Error Analysis for Elliptic Interface Problems on Curved Domains},pdfauthor={Fang Liu, Yuanhui Lin, Tao Lin, Xu Zhang}}
\graphicspath{{assets/}}
\numberwithin{equation}{section}
\newtheorem{theorem}{Theorem}[section]
\newtheorem{lemma}{Lemma}[section]

\newtheorem{remark}{Remark}[section]
\newtheorem{assumption}[theorem]{Assumption}
\newcommand{\Th}{\mathcal{T}_h}
\newcommand{\Eh}{\mathcal{E}_h}
\newcommand{\Eo}{\mathring{\mathcal{E}}_h}
\newcommand{\Es}{\mathcal{E}_h^{s}}
\newcommand{\Hc}{\mathcal{H}}
\newcommand{\KF}{K_F}
\newcommand{\hKF}{\widehat K_F}
\newcommand{\hK}{\widehat K}
\newcommand{\tK}{\widetilde K}
\newcommand{\hG}{\widehat\Gamma_K}
\newcommand{\Pm}{\mathbb{P}_m}

\DeclareMathOperator{\id}{id}
\DeclareMathOperator{\diam}{diam}
\DeclareMathOperator{\dist}{dist}

\newcommand{\commentout}[1]{{}} % for large block comments

\newcommand{\jump}[1]{\left[\!\left[ #1 \right]\!\right]}
\newcommand{\avg}[1]{\{\!\{#1\}\!\}}
\newcommand{\triplenorm}[1]{\left|\!\left|\!\left| #1 \right|\!\right|\!\right|}
\title{A High-Order Geometry-Conforming
Selective Discontinuous Galerkin Method and Its Error Analysis
for Elliptic Interface Problems on Curved Domains}
\author{Fang Liu\thanks{Department of Mathematics, Oklahoma State University, Stillwater, OK 74078, USA. Email: \href{mailto:fliu@okstate.edu}{fliu@okstate.edu}.}
\and Yuanhui Lin\thanks{School of Mathematics, Sichuan University, Chengdu, Sichuan 610065, China. Email: \href{mailto:l1006569522@outlook.com}{l1006569522@outlook.com}.}
\and Tao Lin\thanks{Department of Mathematics, Virginia Tech, Blacksburg, VA 24061, USA. Email: \href{mailto:tlin@vt.edu}{tlin@vt.edu}.}
\and Xu Zhang\thanks{Department of Mathematics, Oklahoma State University, Stillwater, OK 74078, USA. Email: \href{mailto:xzhang@okstate.edu}{xzhang@okstate.edu}.}}
\date{}
\begin{document}
\maketitle

\begin{abstract}
We develop and analyze a high-order geometry-conforming immersed finite element 
(GC-IFE) method for elliptic interface problems on curved domains. The method uses 
a triangular mesh that conforms to the outer boundary but need not fit the internal 
interface. Using Frenet--Serret coordinates, the method combines
GC-IFE spaces on interface elements, geometry-conforming finite element
spaces on curved boundary elements, and standard Lagrange spaces elsewhere. The global degrees of freedom of 
the proposed method are asymptotically the same as that of the standard continuous finite element method on the same mesh. The main analytical contributions are approximation estimates and trace inequalities for triangular GC-IFE functions. These estimates are used to establish stability and 
an optimal $O(h^m)$ error bound in an energy norm for the proposed method, which 
is based on a symmetric selective discontinuous Galerkin formulation. An optimal 
$O(h^{m+1})$ error bound in the $L^2$ norm is also obtained under the standard 
regularity assumption on the dual problem. Numerical experiments 
confirm theoretical convergence rates and demonstrate a reduction 
in of global degrees of freedom to a fully discontinuous
Galerkin method.
\end{abstract}


\section{Introduction}
Elliptic interface problems describe diffusion in heterogeneous media with discontinuous material coefficients. Their solutions typically have reduced global regularity, even when they are smooth within each material subdomain \cite{1963LittmanStampacchiaWeinberger,1998ChenZou}. High-order numerical approximation therefore requires suitable treatment of both the interface jump conditions and the geometry. On a curved domain, resolving the interface alone is insufficient: boundary approximation errors may also limit the overall convergence rate.

\begin{figure}[!htb]
    \centering
    \includegraphics[width=0.4\textwidth]{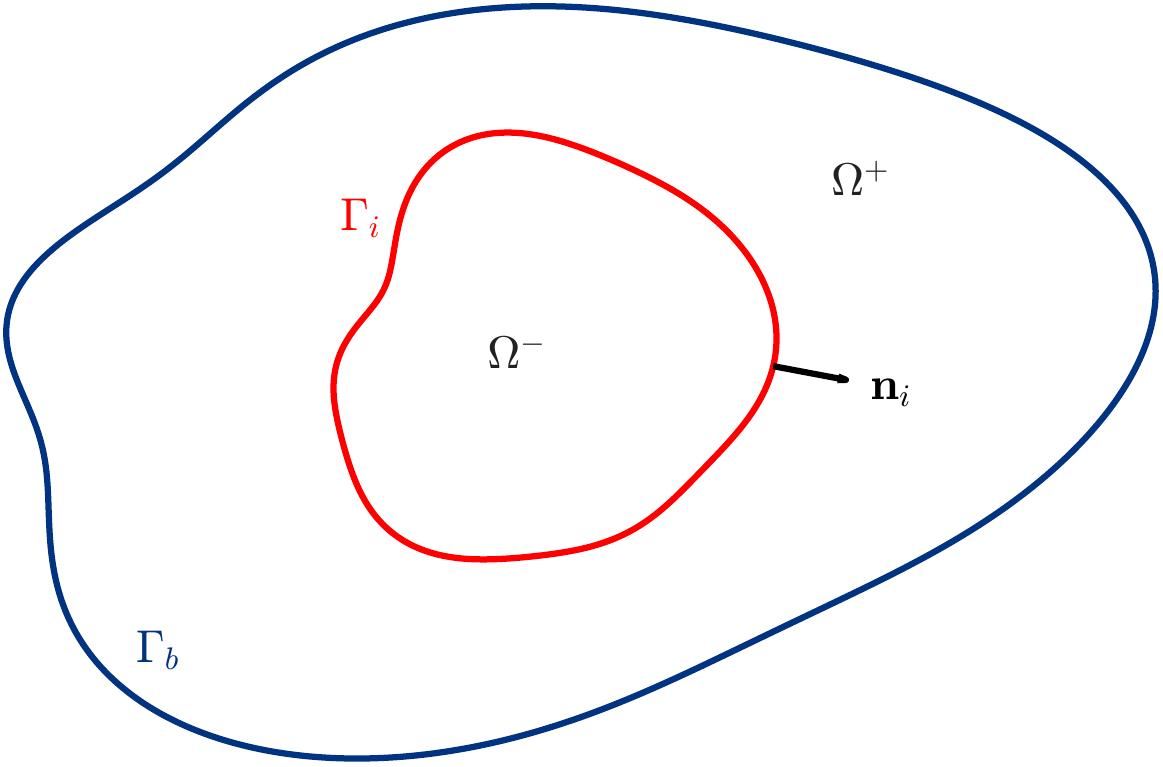}
    \caption{An illustration of the geometry of interface problems: a curved outer boundary $\Gamma_b$ enclosing a curved interface $\Gamma_i$.}
    \label{domain}
\end{figure}

Let $\Omega\subset\mathbb R^2$ be a bounded domain with a curved boundary $\Gamma_b:=\partial\Omega$. Assume that
$\Omega$ is divided into two subdomains $\Omega^-$ and $\Omega^+$ by a smooth closed interface $\Gamma_i \subset \Omega$, as illustrated in Figure~\ref{domain}. We consider the following interface problem:
 \begin{align}
-\nabla \cdot(\beta\nabla u)=f,&\quad\text{ in } \Omega^-\cup \Omega^+,\label{1.1}\\
u=g, &\quad\text{ on } \Gamma_b, \label{1.2}
\end{align}
where $\beta|_{\Omega^\pm}=\beta^\pm>0$ are constants, and $u$ is required to satisfy the following interface jump conditions:
\begin{equation}
  \jump{u}_{\Gamma_i}=0,
  \qquad
  \jump{\beta\partial_{\mathbf n_i}u}_{\Gamma_i}=0.
  \label{1.3}
\end{equation}
Here, $\mathbf n_i$ is the unit normal to $\Gamma_i$ pointing from
$\Omega^-$ to $\Omega^+$, and
$\jump{v}_{\Gamma_i}:=v^+|_{\Gamma_i}-v^-|_{\Gamma_i}$, with
$v^\pm=v|_{\Omega^\pm}$.
For high-order approximations, i.e., $m\geq 2$, we further
assume that the solution satisfies the extended jump conditions
\begin{equation}
\jump{\beta\partial_{\mathbf n_i}^{\,j}\Delta u}_{\Gamma_i}=0,
\qquad 0\le j\le m-2.
  \label{1.4}
\end{equation}
For sufficiently smooth solutions, these conditions follow from
\eqref{1.1} when $f\in C^{m-2}(\Omega)$.

Fitted finite element methods require meshes to align with both the interface and the outer boundary. 
Isoparametric elements are typically used to represent curved element boundaries to high order \cite{1968ErgatoudisIronsZienkiewicz,1986Lenoir}, but generating and updating interface-fitted meshes can be costly, especially for complicated shapes or moving interfaces. Unfitted meshes avoid this constraint, although standard polynomial spaces on such elements generally cannot resolve derivative jumps of the solution, which may compromise accuracy \cite{1970Ivo,2014StefanThomas}.

Immersed finite element (IFE) methods incorporate interface jump conditions into local approximation spaces \cite{2001LinLinRogersRyan,2004LiLinLinRogers}. High-order IFE functions have been developed by least-squares formulations \cite{2017AdjeridGuoLin}, Cauchy extensions \cite{2019GuoLin}, actual-interface constructions \cite{2018AdjeridMohamedLin}, and local enrichments for nonhomogeneous jumps \cite{2023AdjeridIvoGuoLin}. 
More recently, a geometry-conforming IFE (GC-IFE) space was introduced by Adjerid, Lin, and Meghaichi \cite{2024AdjeridLinMeghaichi,2025AdjeridLinMeghaichi,2026AdjeridLinMeghaichi} using local Frenet coordinates to straighten the interface and enforce the physical interface conditions on the prescribed curve. GC-IFE spaces on triangular meshes were
subsequently introduced in \cite{2025LinLinZhang}, where their approximation properties and the convergence of associated discontinuous Galerkin methods were demonstrated numerically.

Accurate approximations on curved domains require suitable finite element spaces 
on boundary elements that contain a portion of the curved outer boundary. In \cite{2026LinLinZhangFeng}, 
the authors developed geometry-conforming finite element (GC-FE) spaces using Frenet transformations 
to straighten the curved edge of a boundary element. That work proposed an IFE method based on standard discontinuous Galerkin (DG) formulations on unfitted meshes, incorporating 
GC-FE spaces on boundary elements, GC-IFE spaces on interface elements, and standard finite element 
spaces on all remaining elements. Although \cite{2026LinLinZhangFeng} numerically demonstrated 
that this IFE method achieves optimal convergence rates, a rigorous error analysis remains to be established.

%Thus, for the unfitted triangular setting, which is the one of practical interest, 
%the approximation and trace estimates needed for a convergence proof have so far been missing.

%Continuous Galerkin (CG) methods share degrees of freedom across neighboring 
%elements and are computationally efficient when local spaces possess compatible traces. 
%For the geometry-conforming spaces considered here, however, independently constructed 
%local functions do not necessarily match along shared element edges. Interior penalty 
%discontinuous Galerkin (DG) formulations accommodate this trace mismatch through weak coupling 
%\cite{2001ArnoldBrezziBernardoMarini,2008Riviere}, but at the cost of introducing 
%separate degrees of freedom for every element. Selective discontinuous Galerkin (SDG) 
%methods mitigate this overhead by retaining continuous coupling wherever possible, 
%applying discontinuous coupling only at selected locations, such as elements near 
%an interface \cite{2014HeLinLin}. A high-order SDG method based on Frenet IFE spaces 
%was recently developed on Cartesian meshes in \cite{2026LiuMeghaichiZhang}. Extending 
%this approach to curved domains requires selective coupling at both boundary and 
%interface elements.

Continuous Galerkin (CG) methods share degrees of freedom across neighboring 
elements and are computationally efficient when local spaces possess compatible traces. 
For the geometry-conforming spaces considered here, however, independently constructed 
local functions do not necessarily match along shared element edges. Interior penalty 
discontinuous Galerkin (DG) methods accommodate this trace mismatch through weak coupling 
\cite{2001ArnoldBrezziBernardoMarini,2008Riviere}, but at the cost of introducing 
separate degrees of freedom for every element. Selective discontinuous Galerkin (SDG) 
methods mitigate this overhead by retaining continuous coupling wherever possible, 
applying discontinuous coupling only at selected locations, such as elements near 
an interface \cite{2014HeLinLin}. A high-order SDG method based on Frenet IFE spaces 
was recently developed on Cartesian meshes in \cite{2026LiuMeghaichiZhang}. However, 
many real-world physical models feature complex, curved outer boundaries. To preserve optimal accuracy for practical applications, it is essential to extend this SDG approach to curved domains, an extension that requires 
selective coupling at both boundary and interface elements.

In this work, we develop and analyze a high-order geometry-conforming selective discontinuous Galerkin 
immersed finite element (GC-SDG-IFE) method for elliptic interface problems on curved 
domains. We employ a mesh that conforms to the outer boundary but remains unfitted 
with respect to the internal interface. We use GC-FE spaces on curved-boundary elements, GC-IFE spaces
on interface-cut elements, and standard Lagrange spaces on regular
elements. These local spaces are coupled via a symmetric SDG formulation: discontinuities are permitted across interior edges adjacent to boundary or interface elements, whereas continuous coupling is retained 
across all remaining interior edges. This formulation unifies the treatment 
of curved outer boundaries and internal interfaces within a single approximation 
framework, reducing the global degrees of freedom compared to a fully 
discontinuous Galerkin discretization.

The primary challenges in the analysis of the proposed SDG-IFE method arise 
on interface elements. Under the inverse Frenet mapping, these elements 
generally yield regions with curved edges, and the internal interface may divide 
an element into subdomains of arbitrary size. Furthermore, because the local basis 
functions are generally non-polynomial in physical coordinates, standard polynomial 
trace inequalities cannot be applied directly. We resolve these issues by introducing 
a comparison between the Frenet elements and corresponding affine triangles, and by 
decomposing each IFE function into piecewise-polynomial components. The main 
contributions of our analysis are summarized as follows:

\begin{itemize}
  \item \textbf{Analysis of triangular GC-IFE spaces.}
  We establish approximation estimates and a discrete flux trace inequality with constants independent of the interface. The proof uses a uniformly sized region on at least one side of the interface.
  % and does not require shape regularity of both sub-elements. 
  These results provide analytical foundations for the triangular GC-IFE spaces developed in \cite{2025LinLinZhang}.
  
  \item \textbf{A GC-SDG-IFE method on curved domains.} This method couples 
   GC-FE functions on boundary elements, GC-IFE functions on interface elements, and standard Lagrange finite element functions on the remaining elements through a symmetric SDG formulation. Only edges adjacent to boundary or interface elements use discontinuous coupling. The remaining interior edges retain CG connectivity. Because discontinuous coupling is confined to boundary and interface elements, the number of global degrees of freedom is asymptotically that of the CG method and smaller than that of a fully DG method.
 
  \item \textbf{Stability and error analysis.}
  We prove the well-posedness of the proposed GC-SDG-IFE method and an
  optimal-order $O(h^m)$ energy error estimate under the stated
  regularity and compatibility assumptions. An $O(h^{m+1})$ estimate
  in the $L^2$ norm follows under the usual regularity assumption for the dual auxiliary problem. The constants in these estimates are explicit in $\beta_{\max} = \max{\beta^\pm}$ and $\beta_{\min}=\min{\beta^\pm}$ and are independent of the location of the interface.
\end{itemize}

%* Organization of this Paper
The rest of this paper is organized as follows. Section \ref{sec2} introduces notation, the Frenet transformations, and the geometric tools for curved Frenet triangles. Section \ref{sec3} develops a high-order geometry-conforming hybrid immersed finite element space and establishes trace inequalities and optimal approximation properties. Section \ref{sec4} introduces a high-order SDG IFE method for curved-domain interface problems. The well-posedness and {\it a priori} error estimates, including optimal-order estimates in the energy norm and the $L^2$ norm, are established. Numerical experiments are presented in Section \ref{sec5} to verify the theoretical results and demonstrate the efficiency and high-order convergence of the proposed method. Finally, brief conclusions are given in Section \ref{sec6}.

%------------------------------------------------------------------
%------------------------------------------------------------------
\section{Preliminaries}\label{sec2}
%------------------------------------------------------------------
%------------------------------------------------------------------

In this section, we recall some notation and the Frenet--Serret apparatus used throughout the paper. Throughout, $m\geq 1$ denotes the polynomial degree of the finite element space and 
\[\beta_{\min}:=\min\{\beta^-,\beta^+\},\qquad \beta_{\max}:=\max\{\beta^-,\beta^+\}.\]
The notation $a\lesssim b$ means $a\leq Cb$, where $C$ is independent of the mesh size, the relative position of the interface, and the diffusion coefficients $\beta^\pm$. We write $a\simeq b$ if $a\lesssim b$ and $b\lesssim a$. For $\iota\in\{+,-\}$, we denote its counterpart by $\iota'$, that is, $\iota'=-$ if $\iota=+$, and $\iota'=+$ if $\iota=-$.

%Let $D\subset\Omega$. For an integer $k\geq 0$ and $1\leq p\leq\infty$, we denote by $W_p^k(\hat{\Omega})$ the standard Sobolev space, equipped with the norm $\|\cdot\|_{W_p^k(\hat{\Omega})}$ and seminorm $|\cdot|_{W_p^k(\hat{\Omega})}$. When $p=2$, we use the notation $H^k(\hat{\Omega})=W_2^k(\hat{\Omega})$. We further denote by $(\cdot,\cdot)_{\partial\hat{\Omega}}$ and $\|\cdot\|_{L^2(\partial\hat{\Omega})}$ the $L^2$ inner product and norm on $\partial\hat{\Omega}$, respectively. 

Let $D\subset\Omega$ be a subdomain. For an integer $k\geq 0$ and $1\leq p\leq\infty$, we denote by $W_p^k(D)$ the standard Sobolev space, equipped with the norm $\|\cdot\|_{W_p^k(D)}$ and seminorm $|\cdot|_{W_p^k(D)}$, and we write $H^k(D)=W_2^k(D)$. We denote by $(\cdot,\cdot)_{\partial D}$ and $\|\cdot\|_{L^2(\partial D)}$ the $L^2$ inner product and norm on $\partial D$. If $D$ intersects the interface $\Gamma_i$, we set $D^{\pm}=D\cap\Omega^{\pm}$ and, for $\alpha>\frac{3}{2}$,
\begin{equation}\label{eq: Hm}
\mathcal{H}^{\alpha}(D)=\Big\{v\in L^2(D) : v|_{D^{\pm}}\in H^{\alpha}(D^\pm),~\jump{v}_{\Gamma_i\cap D}=\jump{\beta\partial_{\textbf{n}_i}v}_{\Gamma_i\cap D}=0\Big\},
\end{equation} 
equipped with the norm $ \|v\|_{\mathcal H^\alpha(D)}^2
 :=\sum_{s=\pm}\|v^s\|_{H^\alpha(D^s)}^2$ 
and the seminorm 
$|v|_{\mathcal H^\alpha(D)}^2
 :=\sum_{s=\pm}|v^s|_{H^\alpha(D^s)}^2$. 
If $D$ does not intersect $\Gamma_i$, we simply set $\mathcal{H}^{\alpha}(D) = H^{\alpha}(D)$. Throughout this paper, $\|\cdot\|$ denotes the Euclidean norm of a vector or the Frobenius norm of a matrix, as appropriate, and for a vector- or matrix-valued function $\mathbf{v}$ we define $\|\mathbf{v}\|_{L^\infty(D)}=\operatorname*{ess\,sup}_{x\in D}\|\mathbf{v}(x)\|$ and $\|\mathbf v\|_{C^0(D)}:=\sup_{\mathbf x\in D}\|\mathbf v(\mathbf x)\|$.
%The same broken-norm convention is used on Frenet rectangles, with the two sides separated by $\eta=0$. We use $\|\cdot\|$ for a Euclidean vector norm or a matrix operator norm.

%equipped with the norm $\|\cdot\|^2_{\mathcal{H}^{\alpha}(\hat{\Omega})}=\|\cdot\|^2_{H^{\alpha}(\hat{\Omega}^+)}+\|\cdot\|^2_{H^{\alpha}(\hat{\Omega}^-)}$, the seminorm
%$|\cdot|^2_{\mathcal{H}^{\alpha}(\hat{\Omega})}=|\cdot|^2_{H^{\alpha}(\hat{\Omega}^+)}+|\cdot|^2_{H^{\alpha}(\hat{\Omega}^-)}$. Throughout this paper, $\|\cdot\|$ denotes the Euclidean norm of a vector or the Frobenius norm of a matrix, as appropriate. For a vector- or
%matrix-valued function $\mathbf{v}$, we define
%$\|\mathbf{v}\|_{L^\infty(\hat{\Omega})}
%=
%\operatorname*{ess\,sup}_{x\in\hat{\Omega}}
%\|\mathbf{v}(x)\|$, and 
%$\|\mathbf v\|_{C^0(\hat{\Omega})}
%:=
%\sup_{\mathbf x\in\hat{\Omega}}
%\|\mathbf v(\mathbf x)\|.$

%------------------------------------------------------------------
\subsection{Mesh Geometry and Frenet Transformations}\label{sec:mesh}
%------------------------------------------------------------------

%We assume that the boundary curve $\Gamma_{\mathrm b}$ and the interface curve $\Gamma_{\mathrm i}$ are disjoint.

Without loss of generality, we assume that the boundary curve $\Gamma_{\mathrm b}$ and the interface curve $\Gamma_{\mathrm i}$ are disjoint. Let $\{\widetilde{\mathcal T}_h\}_{h>0}$ be a family of shape-regular straight-edged triangular meshes of the domain $\Omega$ that is not required to fit the interface $\Gamma_{\mathrm i}$. We assume that $h$ is sufficiently small so that no element intersects both $\Gamma_{\mathrm b}$ and $\Gamma_{\mathrm i}$. We partition $\widetilde{\mathcal T}_h$ into three disjoint subsets, $\widetilde{\mathcal T}_h^b$, $\widetilde{\mathcal T}_h^i$, and $\widetilde{\mathcal T}_h^r$, consisting of boundary, interface, and remaining interior elements, respectively. Specifically, a triangle $\widetilde K$ is called a boundary element if it has an edge whose endpoints lie on $\Gamma_{\mathrm b}$, an interface element if its interior is intersected by $\Gamma_{\mathrm i}$, and a regular element otherwise.

\begin{figure}[H]
    \centering
    
    \begin{subfigure}{0.45\textwidth}
        \centering
       \includegraphics[width=\linewidth]{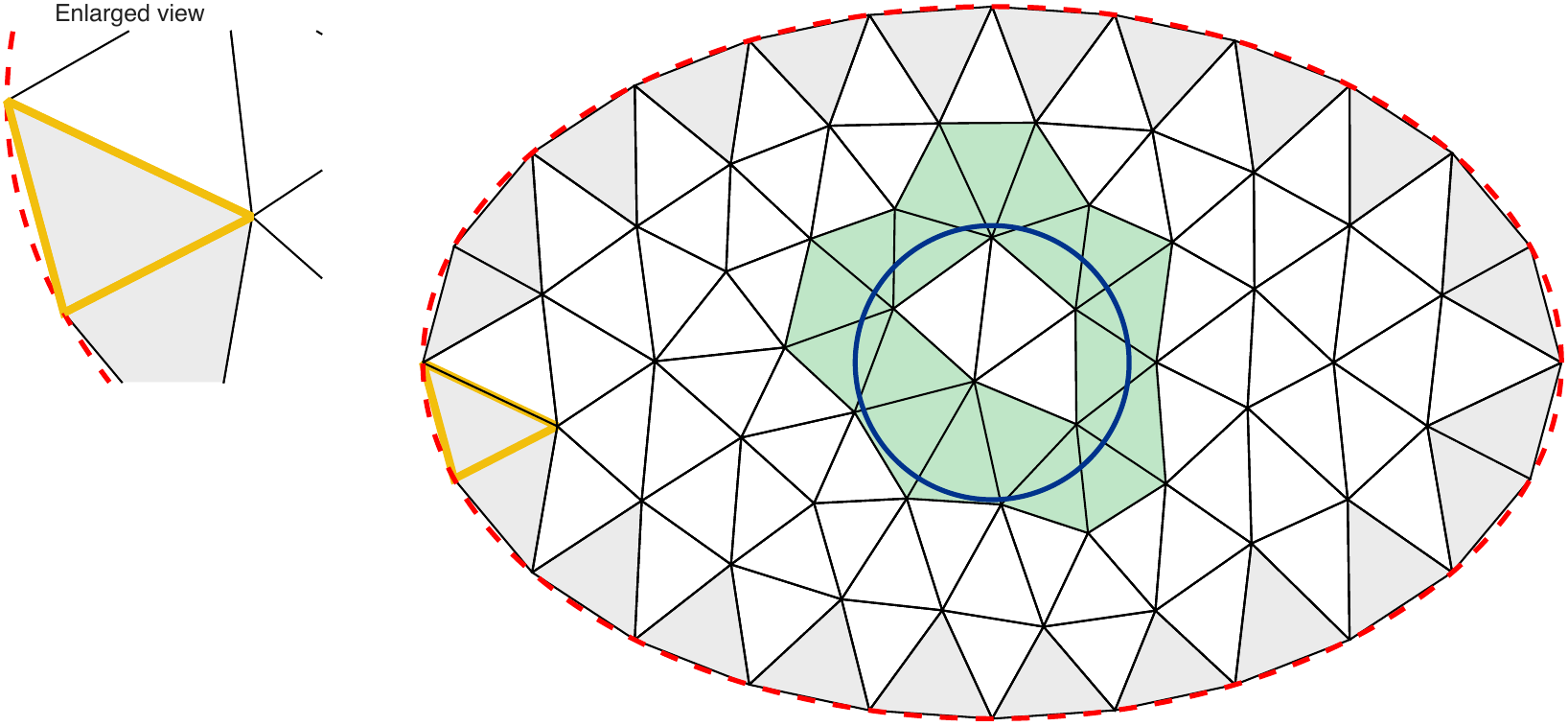}
        \caption{Straight-edged background mesh}
   \end{subfigure}
    \hfill
    \begin{subfigure}{0.45\textwidth}
        \centering
        \includegraphics[width=\linewidth]{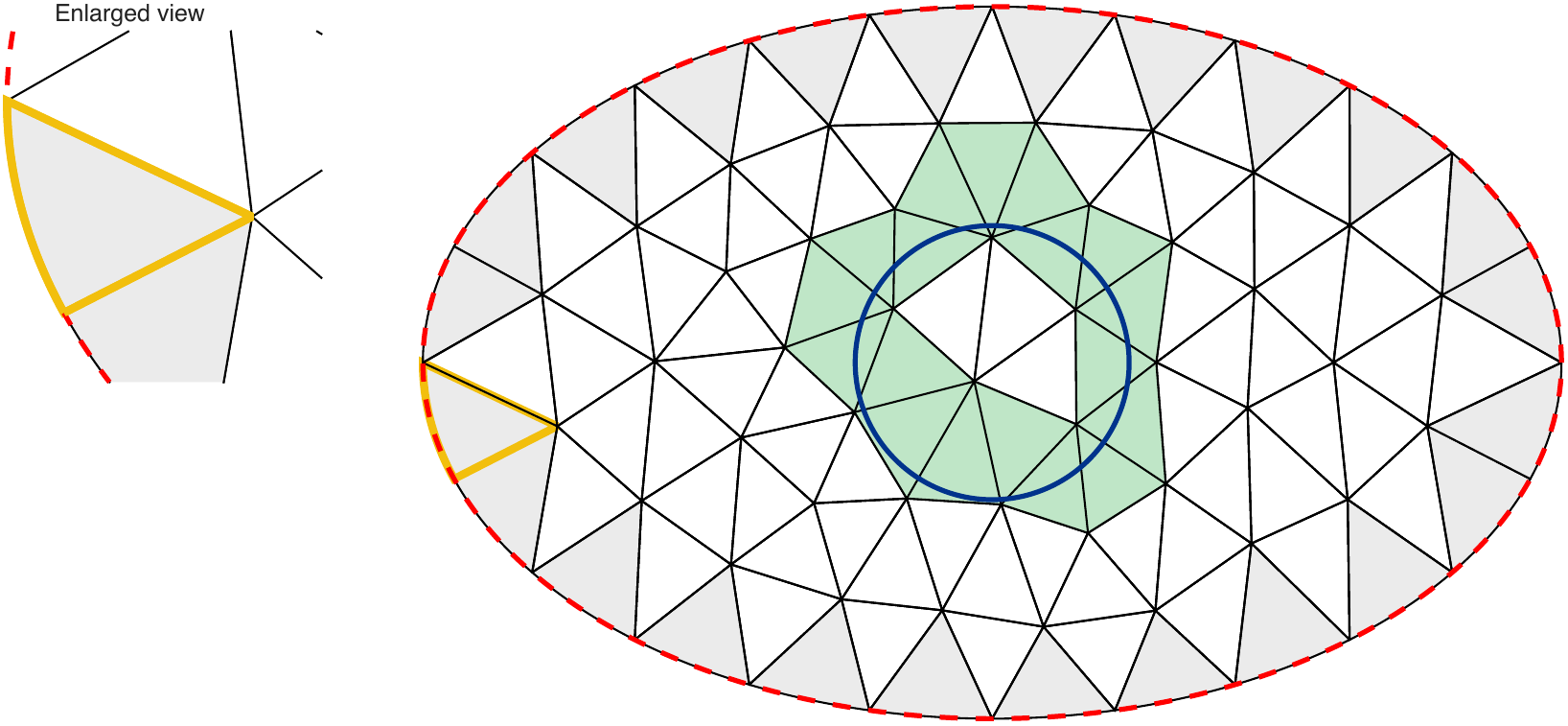}
        \caption{Reconstructed mesh with curved boundary elements}
    \end{subfigure}
    
    \caption{Background mesh $\widetilde{\mathcal{T}}_h$ and the reconstructed mesh $\mathcal{T}_h$. Boundary, interface, and regular elements are shown in grey, green, and white, respectively.} 
    \label{fig:meshdefi}
\end{figure}

For each $\widetilde{K}\in\widetilde{\mathcal{T}}_h^b$, we replace the straight boundary edge with the corresponding arc of $\Gamma_{\mathrm{b}}$. This yields the reconstructed triangulation $\mathcal{T}_h$, whose boundary, interface, and regular element sets are denoted by $\mathcal{T}_h^b$, $\mathcal{T}_h^i$, and $\mathcal{T}_h^r$; thus $\mathcal{T}_h^i=\widetilde{\mathcal{T}}_h^i$, $\mathcal{T}_h^r=\widetilde{\mathcal{T}}_h^r$, whereas each element in $\mathcal{T}_h^b$ is a curved triangle with one curved edge lying on $\Gamma_{\mathrm{b}}$ and two straight edges. Since $\mathcal T_h^i\cap\mathcal T_h^b=\varnothing$, $\mathcal{T}_h=\mathcal{T}_h^i\cup\mathcal{T}_h^b\cup\mathcal{T}_h^r$ is a partition of $\Omega$ (Figure~\ref{fig:meshdefi}). For each $K \in \mathcal{T}_h$, let $h_K := \operatorname{diam}(K)$ and $h:= \max_{K\in\mathcal{T}_h} h_K$. We call $\mathcal T_h^{s}:=\mathcal T_h^b\cup\mathcal T_h^i$ the set of \emph{selected elements}; these are the elements on which non-polynomial local spaces will be used.

For a mesh $\mathcal{T}_h$, let $\mathcal E_h$ be the set of its edges, and let $\mathcal E_h^b$ and $\mathring{\mathcal E}_h$ be its boundary and
interior edge sets, respectively. The set of its \emph{selected edges} is
\[
 \mathcal E_h^s
 :=\{B\in\mathcal E_h:\ B\subset\partial K\text{ for some }
                      K\in\mathcal T_h^s\},
 \qquad
 \mathring{\mathcal E}_h^s:=\mathcal E_h^s\cap\mathring{\mathcal E}_h.
\]
We use $\mathcal E_h^{sb}$ and $\mathring{\mathcal E}_h^{si}$ to denote the sets of edges
of boundary and interface elements, respectively, so
$\mathcal E_h^s=\mathcal E_h^{sb}\cup\mathring{\mathcal E}_h^{si}$.
We also write $\mathring{\mathcal E}_h^i$ for the interior edges cut by
the interface and $\mathring{\mathcal E}_h^n$ for the other interior edges.

For $\alpha>\frac{3}{2}$, we define the following broken immersed Sobolev space:
\begin{equation}\label{C_space}
 \mathcal H^\alpha(\mathcal T_h)
 :=\{v:\ v|_K\in\mathcal H^\alpha(K)\ \forall K\in\mathcal T_h;~~
         \jump v_B=0\ \forall B\in
         {\mathcal E}_h\setminus {\mathcal E}_h^s\}.
\end{equation}

We now recall the Frenet transformation introduced in \cite{2024AdjeridLinMeghaichi}. Let $\Gamma\subset\mathbb{R}^2$ be a regular $C^2$ curve parameterized by $\mathbf{g}:[\xi_s,\xi_e]\to\mathbb{R}^2$, $\mathbf{g}(\xi)=(g_1(\xi),g_2(\xi))^T$. We assume it is also a simple curve without loss of generality. 
%for a closed curve the parametrization is extended periodically with period $\xi_e-\xi_s$. 
The unit tangent, the unit normal, and the signed curvature of $\Gamma$ are
\[\boldsymbol{\tau}(\xi)
=\frac{\mathbf{g}'(\xi)}{\|\mathbf{g}'(\xi)\|},
\qquad
\mathbf{n}(\xi)
=
Q 
\boldsymbol{\tau}(\xi),
\qquad
\kappa(\xi)
=
\frac{\mathbf{g}'(\xi)^T
Q 
\mathbf{g}''(\xi)}
{\|\mathbf{g}'(\xi)\|^3},\qquad
Q  = \begin{bmatrix}0&1\\-1&0\end{bmatrix}.
\]
For sufficiently small $\epsilon>0$, the Frenet transformation
 \[P(\eta,\xi)=\mathbf g(\xi)+\eta\mathbf n(\xi)\]
is one-to-one on $[-\epsilon,\epsilon]\times[\xi_s,\xi_e]$; we denote its image, an $\epsilon$-tubular neighborhood of $\Gamma$ \cite{Abate_Tovena_2012}, by $N(\epsilon)$ and its inverse by $R=P^{-1}$.

When $\Gamma = \Gamma_b$ and $\Gamma = \Gamma_i$, we denote the corresponding Frenet transformations and their inverses by $(P_b,R_b)$ and $(P_i,R_i)$. The parametrization of $\Gamma_{\mathrm b}$ is oriented so that $\mathbf n_b$ is the unit outward normal to $\Omega$, and that of $\Gamma_{\mathrm i}$ so that $\mathbf n_i$ points from $\Omega^-$ into $\Omega^+$. We assume that
\[
h<\min\{\epsilon_b,\epsilon_i\},
\qquad
N_b(h)\cap N_i(h)=\emptyset,
\]
where $\epsilon_b$ and $\epsilon_i$ are the parameters for the tubular neighborhoods of $\Gamma_b$ and $\Gamma_i$, respectively. Then, we restrict the Frenet transformations to the parameter regions
\[
P_b:[-h,0]\times[\xi_{sb},\xi_{eb}]\to N_b(h),
\qquad
P_i:[-h,h]\times[\xi_{si},\xi_{ei}]\to N_i(h),
\]
where $N_b(h)$ and $N_i(h)$ denote the images of the respective parameter regions. In Frenet coordinates, the interface is the line $\eta=0$, and $\eta<0$ ($\eta>0$) corresponds to $\Omega^-$ ($\Omega^+$); accordingly $\hat\beta:=\beta\circ P_i$ equals $\beta^-$ for $\eta<0$ and $\beta^+$ for $\eta>0$.

\begin{figure}[!htb]
    \centering
    \includegraphics[width=0.7\textwidth]{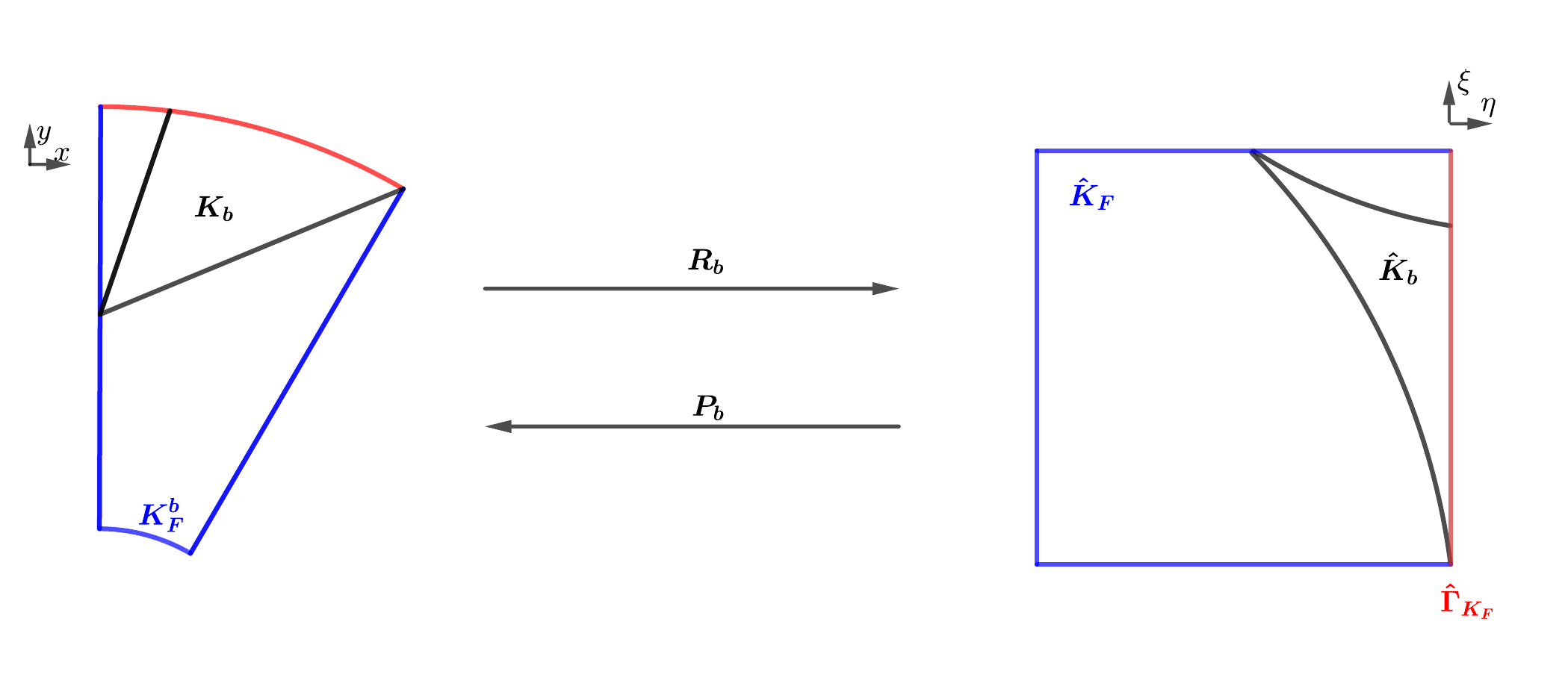}
    \caption{Frenet transformation $P_{b}$ and its inverse $R_{b}$ for a boundary element $K$.}
    \label{fig:frenet_boundary}
\end{figure}

\begin{figure}[!htb]
    \centering
    \includegraphics[width=0.7\textwidth]{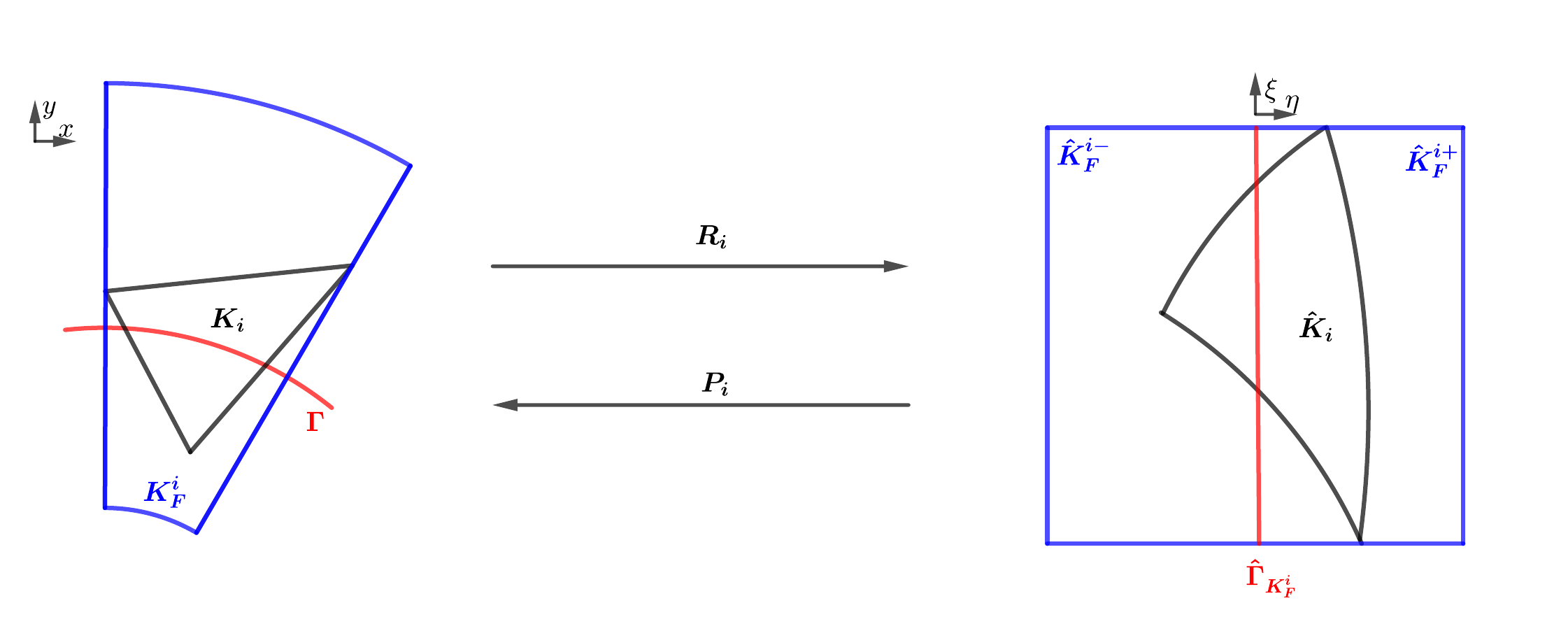}
    \caption{Frenet transformation $P_{i}$ and its inverse $R_{i}$ for an interface element $K$.}
    \label{fig:frenet_interface}
\end{figure}

For $q\in\{b,i\}$ and each $K\in\mathcal{T}_h^q$, let $\xi_k^q$ denote the Frenet parameter of the closest point on $\Gamma_q$ to the vertex $A_k$ of $K$, and let $a_K=\min_{1\le k\le3}\xi_k^q$, $b_K=\max_{1\le k\le3}\xi_k^q$. The \emph{fictitious elements} in Frenet coordinates and in the physical domain are, respectively,
\[\hKF=[-h,0]\times[a_K,b_K]\ (K\in\Th^b),
\qquad
\hKF=[-h,h]\times[a_K,b_K]\ (K\in\Th^i),
\qquad
\KF=P_q(\hKF),\]
and $\hG:=\{0\}\times[a_K,b_K]$ is the segment of $\Gamma_q$ associated with $K$ in Frenet coordinates, so that $|\hG|=b_K-a_K$. We also write $\hK:=R_q(K)\subset\hKF$ for the image of $K$ in Frenet coordinates; for $K\in\Th^i$ it is a \emph{curved} triangle, and $\hK^\pm:=\hK\cap\{\pm\eta>0\}$, $\hKF^{\pm}:=\hKF\cap\{\pm\eta>0\}$. The superscript $q$ is omitted from $\hKF,\KF,\hK$ since it is determined by $K$. Figures~\ref{fig:frenet_boundary} and \ref{fig:frenet_interface} illustrate the Frenet mappings and their inverses.


For the subsequent analysis, let $q\in\{b,i\}$, and assume that $\mathbf{g}_q$ is a regular $C^{m+2}$ parametrization of $\Gamma_q$ satisfying
\begin{equation}\label{eq:geometric_regular}
\|\mathbf{g}_q'(\xi)\|\simeq 1,
\qquad
\|\mathbf{g}_q^{(k)}(\xi)\|\lesssim 1,
\quad 0\le k\le m+2,
\end{equation}
for $\xi\in[\xi_{sq},\xi_{eq}]$. Consequently, $\kappa_q\in C^m([\xi_{sq},\xi_{eq}])$. Furthermore, we assume that
\begin{equation}\label{eq:curvature_condition}
h\kappa_{\max}\le \frac12,
\qquad
\kappa_{\max}
=\max_{q\in\{b,i\}}
\max_{\xi\in[\xi_{sq},\xi_{eq}]}
|\kappa_q(\xi)|.
\end{equation}
Writing $\psi_q(\eta,\xi):=(1+\eta\kappa_q(\xi))^{-1}$, the Frenet formulas $\mathbf n_q'=\kappa_q\|\mathbf g_q'\|\boldsymbol\tau_q$ and $\boldsymbol\tau_q'=-\kappa_q\|\mathbf g_q'\|\mathbf n_q$ give 
\[
\hat D P_q=[\,\mathbf n_q,~~\|\mathbf g_q'\|\psi_q^{-1}\boldsymbol\tau_q\,],
 \qquad
 DR_q=\begin{bmatrix}\mathbf n_q^T\\
                         \psi_q\|\mathbf g_q'\|^{-1}\boldsymbol\tau_q^T\end{bmatrix}
\]
 for the Jacobian matrices of $P_q$ and $R_q$. Under \eqref{eq:geometric_regular}--\eqref{eq:curvature_condition}, $\tfrac12\le\psi_q^{-1}\le\tfrac32$, and therefore
\begin{equation}\label{eq:frenet_jacobian_bound}
|\det \hat{D}P_q|=\|\mathbf g_q'\|\,\psi_q^{-1}\simeq 1,\qquad |\det DR_q|\simeq 1,\qquad \|\hat DP_q\|+\|DR_q\|\lesssim 1,
\qquad q\in\{b,i\},
\end{equation}
on the respective parameter regions; moreover, by \eqref{eq:geometric_regular}, all derivatives of $P_q$ and $R_q$ up to order $m+1$ are uniformly bounded.

Since $K\in\mathcal{T}_h^q$ intersects $\Gamma_q$ and $h_K\le h$, every point of $K$ lies within distance $h$ of $\Gamma_q$; hence $K\subset N_q(h)$ and $|\eta(X)|\le h$ for all $X\in K$. The $\epsilon$-tubular neighborhood property 
\cite{Abate_Tovena_2012} allows us to extend the normal and tangent fields to $K$:
\begin{equation}\label{extend}
\mathbf{n}_q(X)=\mathbf{n}_q(\xi),\quad\boldsymbol{\tau}_q(X)=\boldsymbol{\tau}_q(\xi),\quad \forall X=\mathbf{g}_q(\xi)+\eta\mathbf{n}_q(\xi)\in K.
\end{equation}
Inverting $\hat DP_q$ gives
\[
\nabla \eta(X)=\mathbf{n}_q(\xi),\quad \nabla\xi(X)=\frac{\psi_q(\eta,\xi)\boldsymbol{\tau}_q(\xi)}{\|\mathbf{g}_q'(\xi)\|},\quad \forall X\in K,
\]
and the chain rule yields
\[
\begin{aligned}
\nabla \mathbf{n}_q(X)&=\mathbf{n}_q'(\xi)(\nabla\xi)^T=\kappa_q(\xi)\,\psi_q(\eta,\xi)\,\boldsymbol{\tau}_q(\xi)\boldsymbol{\tau}_q^T(\xi),\\
\nabla \boldsymbol{\tau}_q(X)&=\boldsymbol{\tau}_q'(\xi)(\nabla\xi)^T=-\kappa_q(\xi)\,\psi_q(\eta,\xi)\,\mathbf{n}_q(\xi)\boldsymbol{\tau}_q^T(\xi).
\end{aligned}
\]
Since $\|\boldsymbol\tau_q\boldsymbol\tau_q^T\|=\|\mathbf n_q\boldsymbol\tau_q^T\|=1$ and $|\psi_q|\le 2$ by \eqref{eq:curvature_condition}, we obtain
\begin{equation}\label{Linfty}
\|\nabla \mathbf{n}_q\|_{L^{\infty}(K)}+\|\nabla \boldsymbol{\tau}_q\|_{L^{\infty}(K)}\le 4\kappa_{\max}.
\end{equation}
Finally, $(\eta,\xi)$ are orthogonal coordinates with scale factors $1$ and $\|\mathbf g_q'\|\psi_q^{-1}$, so that the Laplacian takes the form $(\Delta u)\circ P_q=\mathcal L(u\circ P_q)$ with
\begin{equation}\label{eq:LFrenet}
\mathcal L\hat u=\partial_\eta^2\hat u+\kappa_q\psi_q\,\partial_\eta\hat u+\frac{\psi_q}{\|\mathbf g_q'\|}\,\partial_\xi\Big(\frac{\psi_q}{\|\mathbf g_q'\|}\,\partial_\xi\hat u\Big),
\end{equation}
see \cite{2024AdjeridLinMeghaichi}. One property of $\mathcal L$ used below is its structure $\mathcal L=\partial_\eta^2+\kappa_q\psi_q\partial_\eta+\mathcal L_\xi$, where $\mathcal L_\xi$ is defined by the last term in \eqref{eq:LFrenet}, which contains derivatives in $\xi$ only. Likewise, since $\nabla\eta=\mathbf n_q$, the normal derivatives in \eqref{1.3}--\eqref{1.4} become $\eta$-derivatives: $(\partial_{\mathbf n_i}^j w)\circ P_i=\partial_\eta^j(w\circ P_i)$.

The geometric properties of the fictitious elements will be used repeatedly. The following two lemmas summarize results established in \cite{2025AdjeridLinMeghaichi} and \cite{2026LinLinZhangFeng}.

\begin{lemma}\label{lemma1}
For $q\in\{b,i\}$ and any $K\in\mathcal T_h^q$,
\begin{equation}
\left|\hG\right|\simeq h_K,\label{2.4}
\end{equation}
where the equivalence constants are independent of $h$ and of the relative position of $\Gamma_q$.
\end{lemma}

\begin{lemma}\label{boundelem}
There exists an integer $N>0$, independent of $h$, such that
\begin{equation}\label{2.5}
\sum_{q\in\{b,i\}}\sum_{K\in\mathcal{T}_h^q}
\chi_{\KF}(X)
\leq N,
\qquad
\forall\,X\in\Omega,
\end{equation}
where $\chi_{\KF}$ denotes the characteristic function of $\KF$.
\end{lemma}
%------------------------------------------------------------------
\subsection{Trace and Inverse Inequalities on \texorpdfstring{$\hat{K}$}{K-hat}}\label{sec:affine-review}
%------------------------------------------------------------------

\begin{figure}[!htb]
 \centering\includegraphics[width=0.8\textwidth]{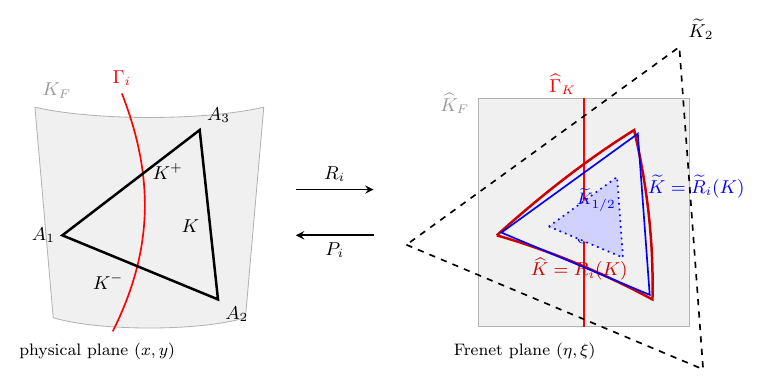}
 \caption{Affine approximation of the Frenet mapping for an interface element. Left: the physical triangle $K$ and its fictitious element $\KF$. Right: in Frenet coordinates the interface is the line $\hG=\{\eta=0\}$, $\hK=R_i(K)$ is a curved triangle (red) and $\tK=\widetilde R_i(K)$ is a straight triangle (blue).}
    \label{fig:sandwich}
\end{figure}

%Inverse and trace inequalities are essential for the coercivity of the bilinear 
%form and for the error estimates. 
In Frenet coordinates, the local spaces on interface elements consist of piecewise polynomials, but the image $\hK=R_i(K)$ of a triangle is a \emph{curved} triangle, so the standard inverse and trace inequalities for polynomials on shape-regular straight elements cannot be invoked directly. In this subsection, we show that $\hK$ is sandwiched between two homothetic copies of a straight, shape-regular triangle, which allows us to transfer polynomial estimates to $\hK$; see Figure~\ref{fig:sandwich}. The analogous analysis for boundary elements can be found in \cite{2026LinLinZhangFeng} and is omitted.

Let $K\in\mathcal T_h^i$, and write $\xi_0=a_K$, $\xi_1=b_K$. Let $\tilde{\mathbf g}_i$ be the linear interpolant of $\mathbf g_i$ at $\xi_0$ and $\xi_1$, and define
\[
\tilde{\boldsymbol{\tau}}_i
=\frac{\tilde{\mathbf g}_i'}{\|\tilde{\mathbf g}'_i\|},
\qquad
\tilde{\mathbf n}_i
=
Q
\tilde{\boldsymbol{\tau}}_i .
\]
The associated affine Frenet mapping and its inverse are $\tilde P_i(\eta,\xi)=\tilde{\mathbf g}_i(\xi)+\eta\,\tilde{\mathbf n}_i$ and $\tilde R_i=\tilde P_i^{-1}$. Let
\[
\tK:=\tilde R_i(K).
\]
By Lemma~\ref{lemma1} and quasi-uniformity, $\xi_1-\xi_0\simeq h$, and hence $\|\tilde{\mathbf g}_i'\|\simeq1$. Since $\tilde{\mathbf g}_i$ is linear, $\tilde{\mathbf g}_i''=\mathbf0$ and, as in \cite{2026LinLinZhangFeng},
\[
|\det(\hat D\tilde P_i)|\simeq|\det(D\tilde R_i)|\simeq1,
\]
with $\hat D\tilde P_i$ and $D\tilde R_i=(\hat D\tilde P_i)^{-1}$ uniformly bounded. Consequently, $\tK$ is the image of the shape-regular triangle $K$ under an affine map with uniformly bounded Jacobian and inverse, hence a shape-regular triangle with $\diam(\tK)\simeq h$, and the standard scaling argument yields, for every $p\in\mathbb P_m(\mathbb R^2)$,
\begin{equation}\label{2.6}
\|\nabla p\|_{L^2(\tK)}
\lesssim h^{-1}\|p\|_{L^2(\tK)},\qquad
\|p\|_{L^2(\partial\tK)}
\lesssim h^{-1/2}\|p\|_{L^2(\tK)},\qquad \|p\|_{L^\infty(\tK)}\lesssim h^{-1}\|p\|_{L^2(\tK)}.
\end{equation}

We introduce a transition map $T=\tilde R_i\circ P_i:\hK\to\tK$ which is close to the identity mapping, as quantified in the following lemma.

\begin{lemma}\label{lemma3}
The following estimates hold for the transition map $T$:
\begin{equation}\label{2.8}
\|T-\id\|_{C^0(\hK)}\lesssim h^2,
\end{equation}
and
\begin{equation}\label{2.9}
\|\hat DT(\hat X)-I_2\|\lesssim h,
\qquad
\left|\det(\hat DT(\hat X))-1\right|\lesssim h,
\qquad \forall \hat X\in\hK,
\end{equation}
where $\mathrm{id}$ is the identity map and $I_2$ is the $2\times 2$ identity matrix.
\end{lemma}

\begin{proof}
Since $\tilde R_i$ is affine, $T-\id=\tilde R_i\circ P_i-\tilde R_i\circ\tilde P_i=D\tilde R_i\,(P_i-\tilde P_i)$, and $P_i-\tilde P_i=(\mathbf g_i-\tilde{\mathbf g}_i)+\eta(\mathbf n_i-\tilde{\mathbf n}_i)$. The linear interpolation error satisfies $\|\mathbf g_i-\tilde{\mathbf g}_i\|_{C^0([\xi_0,\xi_1])}\lesssim h^2$ and $\|\mathbf g_i'-\tilde{\mathbf g}_i'\|_{C^0([\xi_0,\xi_1])}\lesssim h$, whence $\|\mathbf n_i-\tilde{\mathbf n}_i\|\lesssim h$; since $|\eta|\le h$ on $\hK$, \eqref{2.8} follows. Differentiating, $\hat DT-I_2=D\tilde R_i(\hat DP_i-\hat D\tilde P_i)$ with $\hat DP_i-\hat D\tilde P_i=[\,\mathbf n_i-\tilde{\mathbf n}_i\,|\,\mathbf g_i'-\tilde{\mathbf g}_i'+\eta\mathbf n_i'\,]=O(h)$, which gives \eqref{2.9}. See also \cite{2025AdjeridLinMeghaichi} for the rectangular case.
\end{proof}

For $\rho>0$, let $\tK_\rho$ denote the homothetic image of $\tK$ with ratio $\rho$ centered at the centroid $\widetilde G$ of $\tK$; thus $\tK_{1/2}\subset\tK=\tK_1\subset\tK_2$, and each $\tK_\rho$ is a shape-regular triangle of diameter $\simeq h$ for which \eqref{2.6} holds as well. The following lemma, illustrated in Figure~\ref{fig:sandwich}, is the key geometric fact of this subsection. 

%\begin{figure}[!htb]
%    \centering
%    \includegraphics[width=0.4\textwidth]{tri.pdf}
%    \caption{Reference Frenet element $\hat{K}$ (red),  Frenet triangle $\tilde K$ (blue) and  fictitious Frenet triangle $\tilde{\tilde{K}}$ (black). }
%    \label{fig2}
%\end{figure}

\begin{lemma}\label{lem:sandwich}
There exists $h_0>0$, independent of the relative position of the interface, such that for every $K\in\Th^i$ with $h<h_0$,
\[
\tK_{1/2}\subset\hK\subset\tK_{2}.
\]
\end{lemma}

\begin{proof}
Let $C_1$ be the constant in \eqref{2.8}, let $\bar B(\mathbf 0,r)$ be the closed disc of radius $r$ centered at the origin, and let $A+B=\{a+b\mid a\in A,\ b\in B\}$ denote the Minkowski sum for two sets 
$A, B \subset \mathbb R^2$.

Let $G=\tilde P_i(\widetilde G)$. Since $\tilde P_i$ is affine, $G$ is the centroid of $K$, and $\tilde X-\widetilde G=D\tilde R_i\,(X-G)$ for $X\in\partial K$, $\tilde X=\tilde R_i(X)\in\partial\tK$. The uniform bounds on $D\tilde R_i$ and its inverse and the shape regularity of $K$ give
\begin{equation}\label{2.11}
\dist(\widetilde G,\partial\tK)
=\min_{X\in\partial K}\|D\tilde{R}_i\,(G-X)\|
\gtrsim\min_{X\in\partial K}\|G-X\|
\gtrsim h_K\simeq h.
\end{equation}
Each edge of $\tK_\rho$ is parallel to the corresponding edge $e$ of $\tK$ at distance $|\rho-1|\dist(\widetilde G,e)$ from it. Hence, by \eqref{2.11}, there is a constant $c_2>0$, independent of $h$ and of the interface position, such that
\begin{equation}\label{eq:fat-homothety}
\tK+\bar B(\mathbf 0,c_2h)\subset\tK_2,\qquad \tK_{1/2}+\bar B(\mathbf 0,c_2h)\subset\tK .
\end{equation}
Set $h_0:=c_2/C_1$ and let $h<h_0$, so that $C_1h^2<c_2h$.

We now show that $\hK\subset\tK_2$. For $\hat X\in\hK$ we have $T(\hat X)\in\tK$ and $\|\hat X-T(\hat X)\|\le C_1h^2$ by \eqref{2.8}; hence $\hK\subset\tK+\bar B(\mathbf 0,C_1h^2)\subset\tK+\bar B(\mathbf 0,c_2h)\subset\tK_2$ by \eqref{eq:fat-homothety}.

We then show that $\tK_{1/2}\subset\hK$. The map $T$ is a homeomorphism of $\hK$ onto $\tK$, so $\partial\hK=T^{-1}(\partial\tK)$ and, by \eqref{2.8}, $\partial\hK\subset\partial\tK+\bar B(\mathbf 0,C_1h^2)$. By \eqref{eq:fat-homothety}, every point of $\tK_{1/2}$ has distance at least $c_2h>C_1h^2$ from $\partial\tK$; therefore $\tK_{1/2}\cap\partial\hK=\varnothing$. Moreover, $T^{-1}(\widetilde G)$ is an interior point of $\hK$ and $\|T^{-1}(\widetilde G)-\widetilde G\|\le C_1h^2<c_2h$, so $T^{-1}(\widetilde G)\in\tK_{1/2}$. Since $\tK_{1/2}$ is connected, does not meet $\partial\hK$, and contains an interior point of $\hK$, it is contained in the interior of $\hK$.
\end{proof}

We employ Lemma~\ref{lem:sandwich} to transfer polynomial estimates from straight to curved Frenet triangles
in the following lemma.

\begin{lemma}\label{lem:polyK}
Let $h<h_0$ with $h_0$ from Lemma~\ref{lem:sandwich}. For every $K\in\Th^i$ and every $p\in\Pm(\mathbb R^2)$, we have
\begin{align}
&\|p\|_{L^2(\hK)}\simeq\|p\|_{L^2(\tK)}\simeq\|p\|_{L^2(\hKF)},\label{2.12}\\
&\|\nabla p\|_{L^2(\hK)}\lesssim h^{-1}\|p\|_{L^2(\hK)},\label{2.17}\\
&\|p\|_{L^2(\partial\hK)}\lesssim h^{-1/2}\|p\|_{L^2(\hK)}.\label{2.18}
\end{align}
\end{lemma}
\begin{proof}
{Recall that if $D_1,D_2$ are measurable subsets of a disc of radius $Ch$ that both contain a disc of radius $ch$, then $\|p\|_{L^2(D_1)}\simeq\|p\|_{L^2(D_2)}$ for all $p\in\Pm(\mathbb R^2)$, with constants depending only on $m$ and $C/c$; this follows by scaling and the equivalence of norms on the finite-dimensional space $\Pm$ \cite{2008BrennerScott,1978Ciarlet}. The sets $\tK_{1/2},\tK,\tK_2$ and the rectangle $\hKF=[-h,h]\times[a_K,b_K]$ (recall $b_K-a_K\simeq h$ by Lemma~\ref{lemma1}) all satisfy these requirements with constants depending only on the shape regularity. Since $\tK_{1/2}\subset\hK\subset\tK_2$ by Lemma~\ref{lem:sandwich}, $\|p\|_{L^2(\tK_{1/2})}\le\|p\|_{L^2(\hK)}\le\|p\|_{L^2(\tK_2)}$, and \eqref{2.12} follows. 

For \eqref{2.17}, the inverse inequality \eqref{2.6} on $\tK_2$ and \eqref{2.12} give 
\[\|\nabla p\|_{L^2(\hK)}\le\|\nabla p\|_{L^2(\tK_2)}\lesssim h^{-1}\|p\|_{L^2(\tK_2)}\simeq h^{-1}\|p\|_{L^2(\hK)}.\]

For \eqref{2.18}, let $\tilde e$ be an edge of $\tK$ and $\hat e=T^{-1}(\tilde e)$ the corresponding curved edge of $\hK$. By \eqref{2.9}, $\|D(T^{-1})\|\lesssim1$, so the change of variables $\hat e=T^{-1}(\tilde e)$ gives $\|p\|_{L^2(\hat e)}\lesssim\|p\circ T^{-1}\|_{L^2(\tilde e)}$.  For $\tilde Y\in\tilde e$, both $\tilde Y$ and $T^{-1}(\tilde Y)\in\hK$ belong to the convex set $\tK_2$, so the mean value theorem and \eqref{2.8} yield
\begin{equation}\label{mean}
|p(T^{-1}(\tilde Y))-p(\tilde Y)|\le C_1h^2\,\|\nabla p\|_{L^\infty(\tK_2)}\lesssim h^2\cdot h^{-1}\|\nabla p\|_{L^2(\tK_2)}\lesssim\|p\|_{L^2(\tK_2)},
\end{equation}
where we used \eqref{2.6} on $\tK_2$ twice. Consequently, using $|\tilde e|\simeq h$, \eqref{2.6}, \eqref{2.12}, and \eqref{mean}, we have
\[
\begin{aligned}
\|p\|_{L^2(\hat e)}\lesssim\|p\circ T^{-1}\|_{L^2(\tilde e)}
&\lesssim\|p\|_{L^2(\tilde e)}+h^{1/2}\|p\|_{L^2(\tK_2)}\\
&\lesssim h^{-1/2}\|p\|_{L^2(\tK)}+h^{1/2}\|p\|_{L^2(\tK_2)}\lesssim h^{-1/2}\|p\|_{L^2(\hK)}.
\end{aligned}\qedhere
\]

Summing over the three edges proves \eqref{2.18}.}
\end{proof}
\begin{remark}
Lemma~\ref{lem:polyK} is the triangular counterpart of the results obtained for rectangular Frenet elements in \cite{2025AdjeridLinMeghaichi}. The proof above is self-contained and shows that only the two properties \eqref{2.8}--\eqref{2.9} of the transition map and the shape regularity of $\tK$ are needed. Henceforth, we assume $h<h_0$.
\end{remark}

%------------------------------------------------------------------
%------------------------------------------------------------------
\section{Geometry-Conforming Hybrid Immersed Finite Element Space}\label{sec3}
%------------------------------------------------------------------
%------------------------------------------------------------------

In this section, we introduce a geometry-conforming hybrid immersed finite element (GC-HIFE) space and investigate its fundamental properties: its approximation capability and trace inequalities for its local spaces.
%------------------------------------------------------------------
\subsection{Local and Global GC-HIFE Spaces}
%------------------------------------------------------------------
Let $d_m:=\dim\Pm(\mathbb R^2)=\frac{(m+1)(m+2)}{2}$. On a Frenet fictitious element $\hKF$, we use polynomials in the Frenet variables $(\eta,\xi)$; a convenient basis of $\Pm(\mathbb R^2)$, used in \cite{2025LinLinZhang} for the implementation, is
\[
\Big\{\hat q_j\Big(\frac{\eta}{h}\Big)\,\hat p_i\Big(\frac{\xi-\xi_{\rm mid}}{\xi_h}\Big)\ :\ 0\le j\le m,\ 0\le i\le m-j\Big\},
\qquad \xi_{\rm mid}=\dfrac{a_K+b_K}{2},\quad \xi_h=\dfrac{b_K-a_K}{2},
\]
where $\{\hat p_i\}_{i=0}^m$ and $\{\hat q_j\}_{j=0}^m$ are bases for the spaces of univariate polynomials with $\deg\hat p_i=i$, $\deg\hat q_j=j$, $\hat q_1(0)=0$, and $\hat q_j(0)=\hat q_j'(0)=0$ for $2\le j\le m$; the vanishing properties of $\hat q_j$ at $0$ are convenient for enforcing conditions on the interface $\eta=0$. The analysis below does not depend on the choice of basis.

\paragraph{Interface elements.} For $K\in\mathcal{T}_h^i$, the local GC-IFE space $\hat{\mathcal V}_\beta^m(\hKF)$ is the space of piecewise polynomials $\hat\phi$ with $\hat\phi|_{\hKF^{\pm}}\in\mathbb P_m(\hKF^{\pm})$ satisfying
\begin{equation}\label{2.2}
\jump{\hat\phi}_{\hG}=0,\quad
\jump{\hat\beta\,\hat\phi_\eta}_{\hG}=0,
\end{equation}
\begin{equation}\label{2.3}
\int_{\hG}
\jump{
\hat{\beta}\frac{\partial^j}{\partial\eta^j}
\mathcal L(\hat\phi)
}_{\hG}\,\hat v\,ds=0,
\qquad
\forall \hat v\in
\mathbb P_{m-j-2}(\hG),
\quad j=0,1,\ldots,m-2,
\end{equation}
where $\jump{\hat w}_{\hG}=\hat w|_{\hKF^+}-\hat w|_{\hKF^-}$ on $\hG$, and $\hat\beta$ and $\mathcal L$ are as in Section~\ref{sec:mesh}. Conditions \eqref{2.2}--\eqref{2.3} are the Frenet-coordinate form of \eqref{1.3}--\eqref{1.4} on $\hG$: \eqref{2.2} is imposed exactly, whereas \eqref{2.3} is imposed in a weak sense. Specifically, \eqref{2.2} imposes $(m+1)+m$ linear conditions and \eqref{2.3} imposes $\sum_{j=0}^{m-2}(m-j-1)=\frac{m(m-1)}{2}$ linear conditions, i.e., $d_m$ conditions on the $2d_m$ coefficients of $\hat\phi$. It is shown in \cite{2025LinLinZhang} that these conditions are linearly independent, so that $\dim\hat{\mathcal V}_\beta^m(\hKF)=d_m$; in fact, the argument in the proof of Lemma~\ref{lem:decomp} below gives an alternative proof of this fact.

\paragraph{Boundary elements.} Following \cite{2026LinLinZhangFeng}, for $K\in\mathcal T_h^b$ the local GC-FE space on $\hKF$ is simply $\hat{\mathcal W}^{m}(\hKF)=\Pm(\hKF)$, the polynomials in the Frenet variables.

\paragraph{Physical spaces.} We pull back these local spaces to the physical fictitious elements by the inverse Frenet transformations and restrict them to $K$:
\[
\mathcal V_{\beta}^{m}(K)
=
\left\{
\hat v\circ R_i|_K\mid
\hat v\in
\hat{\mathcal V}_{\beta}^{m}(\hKF)
\right\},\quad K\in\Th^i;
\qquad
\mathcal W^{m}(K)
=
\left\{
\hat w\circ R_b|_K\mid
\hat w\in
\Pm(\hKF)
\right\},\quad K\in\Th^b.
\]
{Functions in $\mathcal V_\beta^m(K)$ are piecewise smooth and satisfy \eqref{1.3} \emph{exactly} on the curved interface $\Gamma_i\cap K$, as well as \eqref{1.4} weakly, because $\nabla\eta=\mathbf n_i$; in physical coordinates they are not polynomials.} The local GC-HIFE space on $K\in\mathcal T_h$ is
\begin{equation}\label{eq: GC-HIFE local space}
S_h^m(K)
=
\begin{cases}
\mathbb{P}_m(K),
& K\in\mathcal T_h^r,\\[1mm]
\mathcal V_{\beta}^{m}(K),
& K\in\mathcal T_h^i,\\[1mm]
\mathcal W^{m}(K),
& K\in\mathcal T_h^b,
\end{cases}
\end{equation}
and the global GC-HIFE space is
\begin{equation}\label{eq:Shm}
    S_h^m(\Omega)=\Big\{v\mid~v|_{K}\in S_h^m(K),~\forall K \in \mathcal{T}_h;~v \text{ is continuous across each edge } B\in\Eo\setminus \Es\Big\}.
\end{equation}
By construction, $S_h^m(\Omega)\subset\Hc^{m+1}(\mathcal T_h)$. The fully discontinuous counterpart $\widetilde S_h^m(\Omega)=\{v\mid v|_K\in S_h^m(K),\ \forall K\in\Th\}$ is used only in Remark~\ref{rem:special2} and in the degree-of-freedom comparison of Section~\ref{sec5}.

%------------------------------------------------------------------
\subsection{Approximation Capability}\label{sec:approx}
%------------------------------------------------------------------
We first introduce local projection and interpolation operators on the three types of elements and derive local approximation estimates; the global estimate then follows by summation. The key task is the estimate on interface elements. 

Let $K\in\Th^i$. For $\iota\in\{+,-\}$, let $\hat{\mathcal A}_m^{\iota}:L^2(\hKF^{\iota})\to\mathbb P_m(\hKF^{\iota})$ denote the $L^2$ projection,
\[
(\hat{\mathcal A_m^{\iota}}\hat u^\iota,
\hat v)_{\hKF^{\iota}}
=
(\hat u^\iota,\hat v)_{\hKF^{\iota}},
\qquad
\forall \hat v\in
\mathbb P_m(\hKF^{\iota}),
\]
for $\hat u$ such that $\hat u = \hat u^\pm$ on $\hKF^\pm$, and let $\Theta^{\iota}_m:\mathbb P_m(\hKF^{\iota'})\to\mathbb P_m(\hKF^{\iota})$ be the \emph{IFE extension} operators, i.e., the linear operators such that
\[
(\hat p,\Theta^+_m\hat p)\in\hat{\mathcal V}_\beta^m(\hKF),\quad\forall\,\hat p\in\mathbb P_m(\hKF^{-}),
\qquad
(\Theta^-_m\hat p,\hat p)\in\hat{\mathcal V}_\beta^m(\hKF),\quad\forall\,\hat p\in\mathbb P_m(\hKF^{+}).
\]
These operators are well defined by the unisolvence result of \cite{2025LinLinZhang}. Let $\iota_*\in\{+,-\}$ be such that $\beta^{\iota_*}=\beta_{\min}$ (if $\beta^+=\beta^-$, either choice may be fixed). The local GC-IFE projection $\hat{\mathcal A}_m$ is defined by
\begin{equation}\label{eq:Ahat}
(\hat{\mathcal A}_m\hat u)^{\iota_*}
=
\hat{\mathcal A}_m^{\iota_*}
\hat u^{\iota_*},
\qquad
(\hat{\mathcal A}_m\hat u)^{\iota_*'}
=
\Theta^{\iota_*'}_m
\left(
\hat{\mathcal A}_m^{\iota_*}
\hat u^{\iota_*}
\right),
\end{equation}
so that $\hat{\mathcal A}_m\hat u\in\hat{\mathcal V}_\beta^m(\hKF)$ by construction. In words: we project on the side with the \emph{smaller} coefficient and extend to the other side through the interface conditions; Lemma~\ref{lemma3.3} below shows that the extension amplifies errors at most by the factor $1+\beta_{\min}/\beta_{\max}\le2$, which is what makes the constants in the approximation estimates independent of $\beta^\pm$. On the physical fictitious element, the projection is
\begin{equation}\label{3.4}
\mathcal A_{K,m} u
=
(\hat{\mathcal A}_m\hat u)\circ R_i,
\qquad
\hat u=u\circ P_i .
\end{equation}
For a boundary element $K\in\mathcal T_h^b$, let $\mathcal A^b_{K,m}:L^2(\KF)\to\mathcal W^m(\KF)$ denote the projection onto the local GC-FE space introduced in \cite{2026LinLinZhangFeng}, and for a regular element $K\in\Th^r$ let $I_{K,m}:C^0(K)\to\mathbb{P}_m(K)$ denote the Lagrange interpolation operator. For $u\in\Hc^{m+1}(\Omega)$, we define $\Pi_h^mu$ elementwise by
\begin{equation}\label{Proj}
\Pi_h^mu|_K=
\begin{cases}
\mathcal{A}^b_{K,m}u, & K\in\mathcal{T}_h^b,\\
\mathcal{A}_{K,m}u, & K\in\mathcal{T}_h^i,\\
I_{K,m}u, & K\in\mathcal{T}_h^r.
\end{cases}
\end{equation}
Since every edge in $\Eo\setminus\Es$ is shared by two regular elements, on which the Lagrange interpolants agree, $\Pi_h^m u\in S_h^m(\Omega)$.

For a regular element $K\in\mathcal{T}_h^r$, the standard interpolation estimate gives
\begin{equation}\label{interp}
|u-I_{K,m}u|_{H^j(K)}\lesssim h^{m+1-j}|u|_{H^{m+1}(K)},\qquad 0\le j\le m+1,
\end{equation}
and for a boundary element $K\in\mathcal{T}_h^b$ the corresponding estimate $|u-\mathcal A_{K,m}^bu|_{H^j(\KF)}\lesssim h^{m+1-j}\|u\|_{H^{m+1}(\KF)}$ is proved in \cite[Lemma~3.1]{2026LinLinZhangFeng}.

It remains to establish the approximation property of $\mathcal A_{K,m}$ on an interface element. Since the local space is constructed in Frenet coordinates, the analysis is carried out on $\hKF$ and then transferred to $\KF$. In Frenet coordinates, $\hKF$ is a rectangle of size $2h\times(b_K-a_K)$ with $b_K-a_K\simeq h$, and the interface is the straight segment $\hG$; consequently, the arguments developed for rectangular interface elements in \cite{2024AdjeridLinMeghaichi} apply verbatim to the following three lemmas, and we omit their proofs. Here and below, $\hat u\in\Hc^{m+1}(\hKF)$ is said to satisfy \eqref{1.4} if $\jump{\hat\beta\,\partial_\eta^j\mathcal L\hat u}_{\hG}=0$ for $0\le j\le m-2$, which is the Frenet form of \eqref{1.4} by Section~\ref{sec:mesh}.

\begin{lemma}\label{lemma3.1}\textup{\cite{1978Ciarlet}}
Let $\iota\in \{+,-\}$ and $\hat{u}^\iota\in H^{m+1}(\hKF^{\iota})$. Then
\begin{equation*}
\left|\hat{\mathcal{A}}^{\iota}_m \hat{u}^\iota-\hat{u}^\iota\right|_{H^j(\hKF^{\iota})}\lesssim h^{m+1-j}\left|\hat{u}^\iota\right|_{H^{m+1}( \hKF^{\iota})}, \quad 0\leq j\leq m+1.
\end{equation*}
\end{lemma}

\begin{lemma}\label{lemma3.2}
Let $\hat{u}\in\mathcal{H}^{m+1}(\hKF)$ satisfy \eqref{1.4}. Then
\begin{align*}
&\left\|\hat{\mathcal{A}}_m^{+} \hat u^+-\hat{\mathcal{A}}_m^{-} \hat u^-\right\|_{L^2(\hG)}\lesssim h^{m+\frac{1}{2}}\left|\hat{u}\right|_{\mathcal{H}^{m+1}( \hKF)},\\
&\left\|\beta^+(\hat{\mathcal{A}}_m^{+} \hat u^+)_\eta-\beta^-(\hat{\mathcal{A}}_m^{-} \hat u^-)_{\eta}\right\|_{L^2(\hG)}\lesssim (\beta^++\beta^-)h^{m-\frac{1}{2}}\left|\hat{u}\right|_{\mathcal{H}^{m+1}( \hKF)},\\
&\left\|\beta^+\frac{\partial^j}{\partial \eta^j}\mathcal{L}(\hat{\mathcal{A}}_m^{+} \hat u^+)-\beta^-\frac{\partial^j}{\partial \eta^j}\mathcal{L}(\hat{\mathcal{A}}_m^{-} \hat u^-)\right\|_{L^2(\hG)}\lesssim (\beta^++\beta^-)h^{m-\frac{3}{2}-j}\left|\hat{u}\right|_{\mathcal{H}^{m+1}( \hKF)},\quad 0\le j\le m-2.
\end{align*}
\end{lemma}

\begin{lemma}\label{lemma3.3}
Let $\hat{u}\in\mathcal{H}^{m+1}(\hKF)$ satisfy \eqref{1.4} and let $\iota\in\{+,-\}$. Then
\begin{equation*}
\left|\hat{\mathcal{A}}_m^{\iota} \hat{u}^\iota-\Theta_m^\iota(\hat{\mathcal{A}}_m^{\iota'} \hat{u}^{\iota'})\right|_{H^{j}(\hKF^{\iota})}\lesssim \Big(1+\frac{\beta^{\iota'}}{\beta^{\iota}}\Big)h^{m+1-j}\left|\hat{u}\right|_{\mathcal{H}^{m+1}( \hKF)},\quad 0\leq j\leq m+1.
\end{equation*}
\end{lemma}

Combining Lemmas~\ref{lemma3.1} and \ref{lemma3.3}, we obtain the approximation estimate for $\hat{\mathcal A}_m$. We include the short proof because it shows how the choice of $\iota_*$ in \eqref{eq:Ahat} removes the dependence on $\beta^\pm$.

\begin{lemma}\label{lem:approx-hat}
Let $\hat u\in\mathcal H^{m+1}(\hKF)$ satisfy \eqref{1.4}. Then
\begin{equation}\label{eq:approx-hat}
\left|
\hat{\mathcal A}_m\hat u-\hat u
\right|_{\mathcal H^j(\hKF)}
\lesssim
h^{m+1-j}
|\hat u|_{\mathcal H^{m+1}(\hKF)},
\qquad 0\le j\le m+1.
\end{equation}
\end{lemma}
\begin{proof}
On $\hKF^{\iota_*}$, $\hat{\mathcal A}_m\hat u=\hat{\mathcal A}_m^{\iota_*}\hat u^{\iota_*}$ and Lemma~\ref{lemma3.1} applies directly. On $\hKF^{\iota_*'}$, by \eqref{eq:Ahat}, the triangle inequality, Lemma~\ref{lemma3.1}, and Lemma~\ref{lemma3.3} with $\iota=\iota_*'$,
\begin{align*}
\big|\hat u^{\iota_*'}-\Theta_m^{\iota_*'}(\hat{\mathcal A}_m^{\iota_*}\hat u^{\iota_*})\big|_{H^j(\hKF^{\iota_*'})}
&\le\big|\hat u^{\iota_*'}-\hat{\mathcal A}_m^{\iota_*'}\hat u^{\iota_*'}\big|_{H^j(\hKF^{\iota_*'})}
+\big|\hat{\mathcal A}_m^{\iota_*'}\hat u^{\iota_*'}-\Theta_m^{\iota_*'}(\hat{\mathcal A}_m^{\iota_*}\hat u^{\iota_*})\big|_{H^j(\hKF^{\iota_*'})}\\
&\lesssim\Big(1+\frac{\beta^{\iota_*}}{\beta^{\iota_*'}}\Big)h^{m+1-j}|\hat u|_{\mathcal H^{m+1}(\hKF)}\le 2h^{m+1-j}|\hat u|_{\mathcal H^{m+1}(\hKF)},
\end{align*}
since $\beta^{\iota_*}=\beta_{\min}\le\beta_{\max}=\beta^{\iota_*'}$. Squaring and adding the estimates on the two sides gives \eqref{eq:approx-hat}.
\end{proof}

\begin{lemma}\label{lemma3.5}
Let $K\in\Th^i$ and let $u\in\mathcal H^{m+1}(\KF)$ satisfy \eqref{1.3}--\eqref{1.4}. Then
\begin{equation}
\left|\mathcal{A}_{K,m} u- u\right|_{\mathcal{H}^j(\KF)}\lesssim h^{m+1-j}\left\|u\right\|_{\mathcal{H}^{m+1}( \KF)},\quad 0\leq j\leq m+1.\label{3.5}
\end{equation}
\end{lemma}
\begin{proof}
{Let $\hat u=u\circ P_i$, so that $\mathcal A_{K,m}u-u=(\hat{\mathcal A}_m\hat u-\hat u)\circ R_i$ by \eqref{3.4}. Since the derivatives of $R_i$ and $P_i$ up to order $m+1$ are uniformly bounded and $|\det DR_i|\simeq|\det\hat DP_i|\simeq1$ by \eqref{eq:frenet_jacobian_bound}, the chain rule gives $|v\circ R_i|_{\mathcal H^j(\KF)}\lesssim\|v\|_{\mathcal H^j(\hKF)}$ and $|v\circ P_i|_{\mathcal H^j(\hKF)}\lesssim\|v\|_{\mathcal H^j(\KF)}$ for $0\le j\le m+1$; note that a $j$-th order derivative of a composition involves derivatives of all orders up to $j$, which is why the full norm appears on the right. Hence, by Lemma~\ref{lem:approx-hat} and $h\le1$,
\[
|\mathcal{A}_{K,m} u-u|_{\mathcal{H}^j(\KF)}
\lesssim
\|\hat{\mathcal{A}}_m\hat u-\hat u\|_{\mathcal{H}^j(\hKF)}
\lesssim\sum_{k=0}^{j}h^{m+1-k}|\hat u|_{\mathcal H^{m+1}(\hKF)}
\lesssim
h^{m+1-j}\|u\|_{\mathcal{H}^{m+1}(\KF)}.\qedhere
\]}
\end{proof}

For a piecewise smooth function $v$, we use the broken seminorms $|v|^2_{\Hc^j(\Th)}:=\sum_{K\in\Th}|v|_{\Hc^j(K)}^2$. The following theorem establishes the approximation property of $S_h^m(\Omega)$.
\begin{theorem}\label{thm:approx}
Let $u\in\mathcal{H}^{m+1}(\Omega)$ satisfy \eqref{1.4} and let $\Pi_h^m u\in S_h^m(\Omega)$ be defined by \eqref{Proj}. Then
\begin{equation}
\left| u-\Pi_h^m u\right|_{\mathcal{H}^j(\Th)}\lesssim h^{m+1-j}\|u\|_{\mathcal{H}^{m+1}(\Omega)},\quad 0\leq j \leq m+1.\label{3.6}
\end{equation}
\end{theorem}
\begin{proof}
For $K\in\Th^r$, we use \eqref{interp}. For $K\in\Th^i\cup\Th^b$, we have $K\subset\KF\subset\Omega$, so \eqref{3.5} and \cite[Lemma~3.1]{2026LinLinZhangFeng} give $|u-\Pi^m_hu|_{\Hc^j(K)}\lesssim h^{m+1-j}\|u\|_{\Hc^{m+1}(\KF)}$. Summing the squares over all elements and using the finite overlap of the fictitious elements (Lemma~\ref{boundelem}), $\sum_{K\in\Th^i\cup\Th^b}\|u\|^2_{\Hc^{m+1}(\KF)}\le N\|u\|^2_{\Hc^{m+1}(\Omega)}$, which proves \eqref{3.6}.
\end{proof}

\begin{remark}\label{rem:extended}
The extended conditions \eqref{1.4} enter only through Lemmas~\ref{lemma3.2}--\ref{lemma3.3}: the weak conditions \eqref{2.3} built into the local space are consistent with $u$ only if $u$ satisfies \eqref{1.4}. For $m=1$, the set of conditions \eqref{1.4} is empty and Theorem~\ref{thm:approx} holds for every $u\in\Hc^2(\Omega)$. For $m\ge2$ and solutions violating \eqref{1.4}, optimal order cannot be expected from any IFE space of this type; enriched constructions \cite{2023AdjeridIvoGuoLin} address this situation and are not considered here.
\end{remark}

%------------------------------------------------------------------
\subsection{Trace Inequalities}\label{sec:trace}
%------------------------------------------------------------------
In this subsection, we prove the trace inequality for the local GC-HIFE spaces that will be used for the stability analysis in Section~\ref{sec4}. For interface elements, two difficulties must be overcome: the straightened interface $\hG$ may cut the curved triangle $\hK$ into pieces of arbitrarily different size, and a GC-IFE function is only piecewise polynomial, with a gradient that jumps across $\hG$ by a factor depending on $\beta^\pm$. 

\begin{figure}[!htb]
    \centering
    \includegraphics[width=0.88\textwidth]{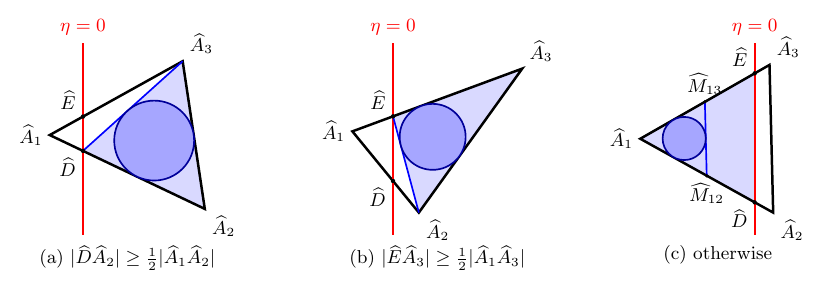}
    \caption{{The three cases in the proof of Lemma~\ref{lem:fat}. The line $\eta = 0$ cuts the triangle at $\widehat{D} \in \widehat{A}_1\widehat{A}_2$ and $\widehat{E} \in \widehat{A}_1\widehat{A}_3$. The shaded region is the sub-triangle used in each case, shown with its incircle.}}
    \label{fig:fatside}
\end{figure}

\begin{lemma}\label{lem:fat}
There exist constants $c_0>0$ and $h_1\in(0,h_0]$, depending only on the shape regularity of the mesh, such that for every $K\in\Th^i$ with $h<h_1$ there is a subelement $\hK^{\iota_K}$ with $\iota_K\in\{+,-\}$ and a disc $B_K$ of radius $c_0h$ such that $B_K\subset\hK^{\iota_K}$. Consequently, for every $p\in\Pm(\mathbb R^2)$ (or $p\in\Pm(\mathbb R^2)^2$),
\begin{equation}\label{eq:fat-equiv}
\|p\|_{L^2(\hKF)}\simeq\|p\|_{L^2(\hK^{\iota_K})}.
\end{equation}
\end{lemma}
\begin{proof}
{Let $\widehat A_1,\widehat A_2$, and $\widehat A_3$ be the vertices of $\hK$ (the images $R_i(A_j)$ of the vertices of $K$), and let $\hK'$ be the straight triangle with these vertices. The vertices of $\tK$ are $T(\widehat A_j)$, so by \eqref{2.8} they lie within $C_1h^2$ of those of $\hK'$. Hence $\hK'$ is shape regular with inradius $r'\simeq h$ for $h$ small, and the curved edge $\hat e=T^{-1}(\tilde e)$ of $\hK$ with endpoints $\widehat A_j,\widehat A_k$ lies within distance $2C_1h^2$ of the straight segment $\widehat A_j\widehat A_k$ (points of $\hat e$ are within $C_1h^2$ of the segment $\tilde e$, whose endpoints are within $C_1h^2$ of $\widehat A_j,\widehat A_k$). Therefore
\begin{equation}\label{eq:fat-boundary}
\partial\hK\subset\partial\hK'+\bar B(\mathbf 0,2C_1h^2).
\end{equation}

\emph{A fat sub-triangle of $\hK'$.} If the line $\{\eta=0\}$ does not meet the interior of $\hK'$, let $\widehat M:=\hK'$ and let $\iota_K$ be the side containing $\hK'$. Otherwise, label the vertices so that $\widehat A_2,\widehat A_3$ lie in one closed half-plane $\{\iota\eta\ge0\}$ and $\widehat A_1$ strictly in the other, and let $\widehat D\in\widehat A_1\widehat A_2$ and $\widehat E\in\widehat A_1\widehat A_3$ be the intersections with $\{\eta=0\}$ (Figure~\ref{fig:fatside}). We distinguish three cases.
\begin{itemize}
\item[(a)] If $|\widehat D\widehat A_2|\ge\frac12|\widehat A_1\widehat A_2|$, let $\widehat M:=\triangle\widehat D\widehat A_2\widehat A_3$ and $\iota_K:=\iota$.
\item[(b)] If $|\widehat E\widehat A_3|\ge\frac12|\widehat A_1\widehat A_3|$, let $\widehat M:=\triangle\widehat E\widehat A_2\widehat A_3$ and $\iota_K:=\iota$.
\item[(c)] Otherwise, $\widehat D$ and $\widehat E$ lie beyond the midpoints $\widehat M_{12},\widehat M_{13}$ of $\widehat A_1\widehat A_2,\widehat A_1\widehat A_3$ as seen from $\widehat A_1$; let $\widehat M:=\triangle\widehat A_1\widehat M_{12}\widehat M_{13}$ and $\iota_K:=\iota'$.
\end{itemize}
In every case, $\widehat M\subset\hK'\cap\{\iota_K\eta\ge0\}$, and the inradius $r_{\widehat M}$ of $\widehat M$ satisfies $r_{\widehat M}\ge\frac12 r'$: in case (c) $\widehat M$ is a homothetic copy of $\hK'$ with ratio $\frac12$, and in cases (a)--(b) $|\widehat M|\ge\frac12|\hK'|$ while the perimeter of $\widehat M$ does not exceed that of $\hK'$, so $r_{\widehat M}=2|\widehat M|/|\partial\widehat M|\ge\frac12 r'$. Let $B'$ be the incircle of $\widehat M$, with center $I$. Since the open disc $B'^{\circ}$ is contained in the interior of $\hK'$ and in the open half-plane $\{\iota_K\eta>0\}$, we have $\dist(I,\partial\hK')\ge r_{\widehat M}$ and $\dist(I,\{\eta=0\})\ge r_{\widehat M}$. Let $B_K$ be the disc centered at $I$ with radius $c_0h:=\frac14 r'\le\frac12 r_{\widehat M}$. Then
\begin{equation}\label{eq:fat-dist}
B_K\subset\{\iota_K\eta>0\},\qquad \dist(B_K,\partial\hK')\ge c_0h .
\end{equation}

\emph{$B_K$ lies in $\hK$.} Let $Z:=\{X\in\hK'\mid \dist(X,\partial\hK')\ge c_0h\}$; this inner parallel set of the convex set $\hK'$ is convex, and $B_K\subset Z$ by \eqref{eq:fat-dist}. Choose $h_1\le h_0$ so small that $2C_1h^2<c_0h$ for $h<h_1$; then $Z\cap\partial\hK=\varnothing$ by \eqref{eq:fat-boundary}. Moreover, the point $T^{-1}(\widetilde G)$ is interior to $\hK$ (Lemma~\ref{lem:sandwich}, Step~3) and lies within $2C_1h^2$ of the centroid of $\hK'$, whose distance to $\partial\hK'$ is $\simeq h$; decreasing $h_1$ if necessary, $T^{-1}(\widetilde G)\in Z$. Being connected, disjoint from $\partial\hK$, and containing an interior point of $\hK$, the set $Z$ lies in the interior of $\hK$. Hence $B_K\subset\hK\cap\{\iota_K\eta>0\}=\hK^{\iota_K}$.

Finally, $B_K\subset\hK^{\iota_K}\subset\hKF$, and both $B_K$ and $\hKF$ are contained in a disc of radius $\lesssim h$ while $B_K$ has radius $c_0h$; the norm-equivalence fact recalled in the proof of Lemma~\ref{lem:polyK} gives $\|p\|_{L^2(\hKF)}\lesssim\|p\|_{L^2(B_K)}\le\|p\|_{L^2(\hK^{\iota_K})}\le\|p\|_{L^2(\hKF)}$.}
\end{proof}

Henceforth, we assume $h<h_1$ so that the results of the preceding lemmas are all valid. 

We now decompose the local GC-IFE space. Let
\[
\begin{aligned}
\mathcal V_0&:=
\Big\{
p\in\Pm(\mathbb R^2):\
p_\eta|_{\hG}=0,\;
\big(\partial_\eta^j\mathcal L p,\hat v\big)_{\hG}=0\
\ \forall \hat v\in\mathbb P_{m-j-2}(\hG),\ 0\le j\le m-2
\Big\},\\
\mathcal Q_m&:=\eta\,\mathbb P_{m-1}(\mathbb R^2),
\end{aligned}
\]
so that $\mathcal Q_m=\{q\in\Pm(\mathbb R^2): q|_{\{\eta=0\}}=0\}$ and $\dim\mathcal Q_m=\frac{m(m+1)}{2}$. A polynomial $p\in\mathcal V_0$ is regarded as a function on $\hKF$ (the same polynomial on both sides), and for $q\in\mathcal Q_m$ we write $\hat\beta^{-1}q$ for the piecewise polynomial equal to $q/\beta^\pm$ on $\hKF^\pm$.

\begin{lemma}\label{lem:decomp}
$\hat{\mathcal V}_\beta^m(\hKF)=\mathcal V_0\oplus\hat\beta^{-1}\mathcal Q_m$, and $\dim\mathcal V_0=m+1$. Consequently, every $\hat\phi\in\hat{\mathcal V}_\beta^m(\hKF)$ can be written uniquely as
\begin{equation}\label{eq:decomp}
\hat\phi=\hat\phi^{(0)}+\hat\beta^{-1}\hat\phi^{(1)},\qquad \hat\phi^{(0)}\in\mathcal V_0,\quad\hat\phi^{(1)}\in\mathcal Q_m,
\qquad\text{and}\qquad \hat\beta\,\hat\nabla\hat\phi=\hat\beta\,\hat\nabla\hat\phi^{(0)}+\hat\nabla\hat\phi^{(1)}.
\end{equation}
\end{lemma}

\begin{proof}
\emph{Inclusions.} Let $p\in\mathcal V_0$. Then $\jump{p}_{\hG}=0$, $\jump{\hat\beta p_\eta}_{\hG}=(\beta^+-\beta^-)p_\eta|_{\hG}=0$, and $\jump{\hat\beta\partial_\eta^j\mathcal Lp}_{\hG}=(\beta^+-\beta^-)\partial_\eta^j\mathcal Lp|_{\hG}$ satisfies \eqref{2.3} by the definition of $\mathcal V_0$; thus $\mathcal V_0\subset\hat{\mathcal V}_\beta^m(\hKF)$. Let $q\in\mathcal Q_m$ and $\hat\phi=\hat\beta^{-1}q$. Then $\jump{\hat\phi}_{\hG}=q|_{\hG}\,(1/\beta^+-1/\beta^-)=0$, and since $\hat\beta$ is constant on each side, $\hat\beta\hat\phi_\eta=q_\eta$ and $\hat\beta\,\partial_\eta^j\mathcal L\hat\phi=\partial_\eta^j\mathcal Lq$ are the restrictions of single smooth functions, whose jumps vanish; thus $\hat\beta^{-1}\mathcal Q_m\subset\hat{\mathcal V}_\beta^m(\hKF)$.

\emph{Trivial intersection.} Let $p\in\mathcal V_0\cap\hat\beta^{-1}\mathcal Q_m$, so $p$ is a polynomial with $p=q/\beta^-$ on $\hKF^-$ and $p=q/\beta^+$ on $\hKF^+$ for some $q\in\mathcal Q_m$. If $\beta^+\ne\beta^-$, then $q/\beta^-=q/\beta^+$ as polynomials, hence $q=0$ and $p=0$. If $\beta^+=\beta^-$, then $p\in\mathcal Q_m\cap\mathcal V_0$, i.e., $p|_{\{\eta=0\}}=0$ and $p_\eta|_{\hG}=0$, so $p=\eta^2s_2$ with $s_2\in\mathbb P_{m-2}$ (and $p=0$ if $m=1$). Suppose $p=\eta^ks_k$ with $s_k\in\mathbb P_{m-k}$ for some $2\le k\le m$. By the structure \eqref{eq:LFrenet} of $\mathcal L$, $\mathcal L(\eta^ks_k)=k(k-1)\eta^{k-2}s_k+\eta^{k-1}r$ with $r$ smooth, so $\partial_\eta^{k-2}\mathcal Lp|_{\hG}=k!\,s_k(0,\cdot)$. The condition in $\mathcal V_0$ with $j=k-2$ and $\hat v=s_k(0,\cdot)\in\mathbb P_{m-k}(\hG)$ gives $s_k(0,\cdot)=0$ on $\hG$, hence $s_k=\eta s_{k+1}$ with $s_{k+1}\in\mathbb P_{m-k-1}$. By induction, $p=\eta^{m+1}s_{m+1}$, which forces $p=0$ since $\deg p\le m$.

\emph{Dimensions.} $\mathcal V_0$ is the kernel of $m+\sum_{j=0}^{m-2}(m-j-1)=\frac{m(m+1)}{2}$ linear functionals on $\Pm(\mathbb R^2)$, so $\dim\mathcal V_0\ge d_m-\frac{m(m+1)}{2}=m+1$ and $\dim(\mathcal V_0\oplus\hat\beta^{-1}\mathcal Q_m)\ge d_m$. On the other hand, the map $\hat\phi\mapsto\hat\phi^-$ is injective on $\hat{\mathcal V}_\beta^m(\hKF)$: if $\hat\phi^-=0$, then $p:=\hat\phi^+$ vanishes on $\hG$ together with $p_\eta$ and satisfies $(\partial_\eta^j\mathcal Lp,\hat v)_{\hG}=0$ for all $\hat v\in\mathbb P_{m-j-2}(\hG)$, so $p=0$ by the induction just performed. Hence $\dim\hat{\mathcal V}_\beta^m(\hKF)\le d_m$, which forces $\hat{\mathcal V}_\beta^m(\hKF)=\mathcal V_0\oplus\hat\beta^{-1}\mathcal Q_m$ and $\dim\mathcal V_0=m+1$. The identity for $\hat\beta\hat\nabla\hat\phi$ in \eqref{eq:decomp} holds because $\hat\beta$ is piecewise constant.
\end{proof}

\begin{remark}
The injectivity argument in the proof also shows that the $d_m$ conditions \eqref{2.2}--\eqref{2.3} are linearly independent, i.e., it reproves the unisolvence of the local GC-IFE space established in \cite{2025LinLinZhang}; it only uses the structure $\mathcal L=\partial_\eta^2+\kappa_i\psi_i\partial_\eta+\mathcal L_\xi$ of the Laplacian in Frenet coordinates.
\end{remark}

The decomposition is stable in the $H^1$ seminorm on each half of the Frenet rectangle. Since $\hKF^{\iota}$ is a rectangle with side lengths $h$ and $b_K-a_K\simeq h$, the following statement coincides with the one proved for rectangular Frenet elements in \cite{2025AdjeridLinMeghaichi}, and we refer to that paper for the proof (which also handles the dependence of $\mathcal V_0$ on $h$ through the coefficients of $\mathcal L$).

\begin{lemma}\label{lem:stab-rect}
Let $\iota\in\{+,-\}$. Then
\begin{equation*}
\|\hat\nabla\hat\phi^{(0)}\|_{L^2(\hKF^{\iota})}
+
\|\hat\nabla\hat\phi^{(1)}\|_{L^2(\hKF^{\iota})}
\lesssim
\|\hat\nabla\hat\phi^{(0)}
+\hat\nabla\hat\phi^{(1)}\|_{L^2(\hKF^{\iota})},
\qquad\forall \hat\phi^{(0)}\in\mathcal V_0,\ \hat\phi^{(1)}\in\mathcal Q_m .
\end{equation*}
\end{lemma}

{We now transfer this stability to the curved triangle $\hK$, with explicit dependence on $\beta^\pm$. This is the key step for the trace inequality.}

\begin{lemma}\label{lem:stab-K}
For every $\hat\phi\in\hat{\mathcal V}_\beta^m(\hKF)$ with the decomposition \eqref{eq:decomp},
\begin{equation}\label{eq:stab-K}
\beta_{\max}\|\hat\nabla\hat\phi^{(0)}\|_{L^2(\hK)}
+
\|\hat\nabla\hat\phi^{(1)}\|_{L^2(\hK)}
\lesssim
\frac{\beta_{\max}}{\sqrt{\beta_{\min}}}
\left\|\sqrt{\hat\beta}\,\hat\nabla\hat\phi\right\|_{L^2(\hK)}.
\end{equation}
\end{lemma}
\begin{proof}
{Let $\iota:=\iota_K$ be the side provided by Lemma~\ref{lem:fat}. On $\hKF^{\iota}$, we have $\hat\phi=\hat\phi^{(0)}+(\beta^{\iota})^{-1}\hat\phi^{(1)}$ with $(\beta^\iota)^{-1}\hat\phi^{(1)}\in\mathcal Q_m$, so Lemma~\ref{lem:stab-rect} gives
\[
\|\hat\nabla\hat\phi^{(0)}\|_{L^2(\hKF^{\iota})}+(\beta^\iota)^{-1}\|\hat\nabla\hat\phi^{(1)}\|_{L^2(\hKF^{\iota})}\lesssim\|\hat\nabla\hat\phi\|_{L^2(\hKF^{\iota})}.
\]
All three gradients are polynomial vectors on $\hKF^\iota$. Since $\hKF^\iota$ is a rectangle with sides $\simeq h$, \eqref{2.12} and the norm-equivalence fact of Lemma~\ref{lem:polyK} give $\|\hat\nabla\hat\phi^{(k)}\|_{L^2(\hK)}\simeq\|\hat\nabla\hat\phi^{(k)}\|_{L^2(\hKF^{\iota})}$, $k=0,1$, while \eqref{eq:fat-equiv} gives $\|\hat\nabla\hat\phi\|_{L^2(\hKF^\iota)}\le\|\hat\nabla\hat\phi^{\iota}\|_{L^2(\hKF)}\lesssim\|\hat\nabla\hat\phi\|_{L^2(\hK^{\iota})}$, where $\hat\phi^\iota\in\Pm$ denotes the polynomial coinciding with $\hat\phi$ on $\hKF^\iota$. Therefore
\[
\begin{aligned}
\beta_{\max}\|\hat\nabla\hat\phi^{(0)}\|_{L^2(\hK)}+\|\hat\nabla\hat\phi^{(1)}\|_{L^2(\hK)}
&\lesssim(\beta_{\max}+\beta^\iota)\|\hat\nabla\hat\phi\|_{L^2(\hK^\iota)}\\
&\le\frac{2\beta_{\max}}{\sqrt{\beta^\iota}}\big\|\sqrt{\hat\beta}\hat\nabla\hat\phi\big\|_{L^2(\hK^\iota)}
\le\frac{2\beta_{\max}}{\sqrt{\beta_{\min}}}\big\|\sqrt{\hat\beta}\hat\nabla\hat\phi\big\|_{L^2(\hK)}.
\end{aligned}\qedhere
\]}
\end{proof}

\begin{remark}
The factor ${\beta_{\max}}/{\sqrt{\beta_{\min}}}$ in \eqref{eq:stab-K} cannot be removed in general. For $m=1$ and an element $K$ lying almost entirely on the side with the larger coefficient $\beta_{\max}$, the IFE function $\hat\phi=\eta$ on that side has $\hat\phi_\eta=\beta_{\max}/\beta_{\min}$ on the other side, so an estimate of $\|\sqrt{\hat\beta}\hat\nabla\hat\phi\|_{L^2(\hKF)}$ by $\|\sqrt{\hat\beta}\hat\nabla\hat\phi\|_{L^2(\hK)}$ necessarily involves the contrast $\beta_{\max}/\beta_{\min}$. Lemma~\ref{lem:stab-K} avoids any estimate on the fictitious element and works on $\hK$ directly.
\end{remark}

Combining the trace estimates on regular, interface, and boundary elements, we obtain the trace inequality for the local GC-HIFE spaces.

\begin{theorem}\label{thm:trace}
For any element $K\in\mathcal T_h$,
\begin{equation}\label{3.10}
\|\beta\nabla v\|_{L^2(\partial K)}
\lesssim
\frac{\beta_{\max}}{\sqrt{\beta_{\min}}}
h^{-1/2}
\|\sqrt{\beta}\nabla v\|_{L^2(K)},
\qquad
\forall v\in S_h^m(K).
\end{equation}
\end{theorem}
\begin{proof}
{If $K\in\Th^r$, then $\beta\equiv\beta^\iota$ on $K$ and the standard polynomial trace inequality \cite{2008BrennerScott} gives $\|\beta\nabla v\|_{L^2(\partial K)}\lesssim\beta^\iota h^{-1/2}\|\nabla v\|_{L^2(K)}=\sqrt{\beta^\iota}\,h^{-1/2}\|\sqrt\beta\nabla v\|_{L^2(K)}$, and $\sqrt{\beta^\iota}\le\sqrt{\beta_{\max}}\le\beta_{\max}/\sqrt{\beta_{\min}}$. If $K\in\Th^b$, then $S_h^m(K)=\mathcal W^m(K)$, $\beta$ is constant on $K$, and the estimate is proved in \cite{2026LinLinZhangFeng}.

Let $K\in\Th^i$ and $v=\hat\phi\circ R_i|_K$ with $\hat\phi\in\hat{\mathcal V}_\beta^m(\hKF)$ decomposed as in \eqref{eq:decomp}. On the curved triangle $\hK$, the identity $\hat\beta\hat\nabla\hat\phi=\hat\beta\hat\nabla\hat\phi^{(0)}+\hat\nabla\hat\phi^{(1)}$, the polynomial trace inequality \eqref{2.18} applied to the polynomial vectors $\hat\nabla\hat\phi^{(0)},\hat\nabla\hat\phi^{(1)}$, and Lemma~\ref{lem:stab-K} yield
\begin{align*}
\|\hat\beta\hat\nabla\hat\phi\|_{L^2(\partial\hK)}
&\leq \beta_{\max}\|\hat\nabla\hat\phi^{(0)}\|_{L^2(\partial\hK)}+\|\hat\nabla\hat\phi^{(1)}\|_{L^2(\partial\hK)}\\
&\lesssim h^{-\frac{1}{2}}\Big(\beta_{\max}\|\hat\nabla\hat\phi^{(0)}\|_{L^2(\hK)}+\|\hat\nabla\hat\phi^{(1)}\|_{L^2(\hK)}\Big)
\lesssim\frac{\beta_{\max}}{\sqrt{\beta_{\min}}}h^{-\frac{1}{2}}\big\|\sqrt{\hat\beta}\hat\nabla\hat\phi\big\|_{L^2(\hK)}.
\end{align*}
Finally, $\nabla v=(DR_i)^T(\hat\nabla\hat\phi)\circ R_i$ and $\hat\nabla\hat\phi=(\hat DP_i)^T(\nabla v)\circ P_i$, and by \eqref{eq:frenet_jacobian_bound} the Jacobian matrices and the area and arc-length elements of $R_i$ and $P_i$ are uniformly bounded; hence $\|\beta\nabla v\|_{L^2(\partial K)}\lesssim\|\hat\beta\hat\nabla\hat\phi\|_{L^2(\partial\hK)}$ and $\|\sqrt{\hat\beta}\hat\nabla\hat\phi\|_{L^2(\hK)}\lesssim\|\sqrt\beta\nabla v\|_{L^2(K)}$, which proves \eqref{3.10}.}
\end{proof}

The error analysis also requires a trace inequality for the exact solution on interface elements, where $u$ is only piecewise $H^2$ but satisfies the jump conditions \eqref{1.3}. Its proof uses the fact that the vector $(\beta\partial_{\mathbf n_i}u,\beta_{\min}\partial_{\boldsymbol\tau_i}u)$ is continuous across the interface, together with the bounds \eqref{Linfty} on the extended normal and tangent fields.

\begin{theorem}\label{thm:traceH2}
Let $K\in\mathcal T_h^i$ and $u\in\mathcal H^2(K)$. Then, for every edge $B$ of $K$,
\begin{equation}\label{eq:trace_beta_grad_u}
\|\beta\nabla u\|_{L^2(B)}
\lesssim
\frac{\beta_{\max}}{\beta_{\min}}
\left(
h^{-1/2}\|\beta\nabla u\|_{L^2(K)}
+h^{1/2}\|\beta D^2u\|_{L^2(K^-\cup K^+)}
\right),
\end{equation}
where $\|\beta D^2u\|_{L^2(K^-\cup K^+)}^2:=\sum_{\iota=\pm}(\beta^\iota)^2\|D^2u\|_{L^2(K^\iota)}^2$.
\end{theorem}
\begin{proof}
Let $\mathbf n=\mathbf n_i$ and $\boldsymbol\tau=\boldsymbol\tau_i$ be the fields extended to $K$ by \eqref{extend}, and define $v_n:=\beta\nabla u\cdot\mathbf n$, $v_t:=\beta_{\min}\nabla u\cdot\boldsymbol\tau$. On $\Gamma_i\cap K$, the flux condition gives $\jump{v_n}=0$, and $\jump{u}=0$ implies that the tangential derivatives of $u^\pm$ coincide, so $\jump{v_t}=0$; since $u\in\Hc^2(K)$, we conclude $v_n,v_t\in H^1(K)$. Let $\mathbf v:=(v_n,v_t)$. Then $|\mathbf v|^2=\beta^2|\nabla u\cdot\mathbf n|^2+\beta_{\min}^2|\nabla u\cdot\boldsymbol\tau|^2$, so
\[
\beta_{\min}^2|\nabla u|^2\leq |\mathbf v|^2\leq \beta^2|\nabla u|^2,
\]
and hence
\[
\|\beta\nabla u\|_{L^2(B)}\leq\beta_{\max}\|\nabla u\|_{L^2(B)}\le \frac{\beta_{\max}}{\beta_{\min}}\|\mathbf v\|_{L^2(B)},\qquad \|\mathbf v\|_{L^2(K)}\le\|\beta\nabla u\|_{L^2(K)}.
\]
By the product rule,
$\nabla v_n=\beta D^2u\,\mathbf n+\beta(\nabla\mathbf n)^T\nabla u$ and $\nabla v_t=\beta_{\min} D^2u\,\boldsymbol\tau+\beta_{\min} (\nabla\boldsymbol\tau)^T\nabla u$ on $K^\pm$, so \eqref{Linfty} gives
\[
\|\nabla\mathbf v\|_{L^2(K)}
\lesssim
\|\beta D^2u\|_{L^2(K^-\cup K^+)}
+
\kappa_{\max}\|\beta\nabla u\|_{L^2(K)}.
\]
Applying the standard trace inequality $\|\mathbf v\|_{L^2(B)}\lesssim h^{-1/2}\|\mathbf v\|_{L^2(K)}+h^{1/2}\|\nabla\mathbf v\|_{L^2(K)}$ for $\mathbf v\in H^1(K)^2$ on the straight triangle $K$ and using $h\kappa_{\max}\le\frac12$ from \eqref{eq:curvature_condition}, we obtain
\begin{align*}
\|\beta\nabla u\|_{L^2(B)}
&\lesssim
\frac{\beta_{\max}}{\beta_{\min}}
\left(
(1+h\kappa_{\max})h^{-1/2}\|\beta\nabla u\|_{L^2(K)}
+
h^{1/2}\|\beta D^2u\|_{L^2(K^-\cup K^+)}
\right)\\
&\lesssim\frac{\beta_{\max}}{\beta_{\min}}
\left(
h^{-1/2}\|\beta\nabla u\|_{L^2(K)}
+
h^{1/2}\|\beta D^2u\|_{L^2(K^-\cup K^+)}\right).\qedhere
\end{align*}
\end{proof}

 %--------------------------------------------------------------------------
%--------------------------------------------------------------------------
\section{The Geometry-Conforming SDG Method}\label{sec4}
%--------------------------------------------------------------------------
%--------------------------------------------------------------------------

In this section, we develop a symmetric selective discontinuous Galerkin (SDG) method based on the GC-HIFE space $S_h^m(\Omega)$, prove its well-posedness, and derive \emph{a priori} error estimates in the energy and $L^2$ norms.

%--------------------------------------------------------------------------
\subsection{The Symmetric SDG Formulation}\label{sec:sdg}
%--------------------------------------------------------------------------
For every interior edge $B\in\Eo$, let $K_{B,1}$ and $K_{B,2}$ denote the two elements sharing $B$, and let $\mathbf n_B$ be the unit normal to $B$ pointing from $K_{B,1}$ into $K_{B,2}$. For a function $u$ defined on $K_{B,1}\cup K_{B,2}$, the average and jump of $u$ on $B$ are
\[
\avg{u}_B = \frac{1}{2} \Big( (u|_{K_{B,1}})|_B + (u|_{K_{B,2}})|_B \Big),
\qquad
\jump{u}_B = (u|_{K_{B,1}})|_B - (u|_{K_{B,2}})|_B .
\]
For $B\in\mathcal{E}_h^b$, $\mathbf n_B$ is the outward unit normal and $\avg{u}_B =\jump{u}_B = (u|_{K_{B,1}})|_B$. {Let $\sigma>0$ be a fixed constant, to be chosen sufficiently large, as specified in Lemma~\ref{lem:coer}, and set
\begin{equation}\label{eq:sigma0}
\sigma_0:=\sigma\,\frac{\beta_{\max}^2}{\beta_{\min}}.
\end{equation}
On the broken space $V_h:=\Hc^{2}(\Th)$, which contains $\Hc^{m+1}(\Th)$ and $S_h^m(\Omega)$ with $m \geq 1$, we define the symmetric bilinear form}
\begin{equation}
\begin{aligned}
a_h(u,v):=
 &\sum_{K\in\mathcal{T}_h}\int_{K}\beta\nabla u\cdot \nabla v\,dX-\sum_{B\in\Es}\int_B\avg{\beta\nabla u\cdot\mathbf{n}_B}\jump{v}\,ds\\
&-\sum_{B\in\Es}\int_B\avg{\beta\nabla v\cdot\mathbf{n}_B}\jump{u}\,ds+\sum_{B\in\Es}\int_B\frac{\sigma_0}{h_B}\jump{u}\jump{v}\,ds,\label{4.2}
\end{aligned}
\end{equation}
and the linear form
\begin{equation}
L_h(v):=\sum_{K\in\mathcal{T}_h}\int_Kfv\,dX+\sum_{B\in\mathcal{E}_h^b}\int_{B}\Big(-\beta\nabla v\cdot\mathbf{n}_B+\frac{\sigma_0}{h_B}v\Big)g\,ds.\label{4.3}
\end{equation}
{Let $u$ be the solution of \eqref{1.1}--\eqref{1.3} and assume $u\in\Hc^2(\Omega)$. Integrating by parts on $K^\pm$ for every element and using that $\beta\nabla u\cdot\mathbf n_B$ is continuous across all interior edges, that $\jump{u}=0$ on $\Gamma_i$ and across all interior edges, and that $\jump{v}=0$ on $\Eo\setminus\Es$ for $v\in V_h$, one obtains the consistency identity}
\begin{equation}
a_h(u,v)=L_h(v),\quad\forall v\in V_h. \label{4.1}
\end{equation}
The symmetric SDG method reads: find $u_h\in S_h^m(\Omega)$ such that
\begin{equation}\label{SDG}
a_h(u_h,v_h)=L_h(v_h),\quad\forall v_h\in S_h^m(\Omega).
\end{equation}
Since $S_h^m(\Omega)\subset V_h$, \eqref{4.1} and \eqref{SDG} imply the Galerkin orthogonality
\begin{equation}\label{eq:GO}
a_h(u-u_h,v_h)=0,\qquad\forall v_h\in S_h^m(\Omega).
\end{equation}

\begin{figure}[!htb]
    \centering

    % ---------- First row ----------
    \begin{subfigure}{0.3\textwidth}
        \centering
        \includegraphics[width=\linewidth]{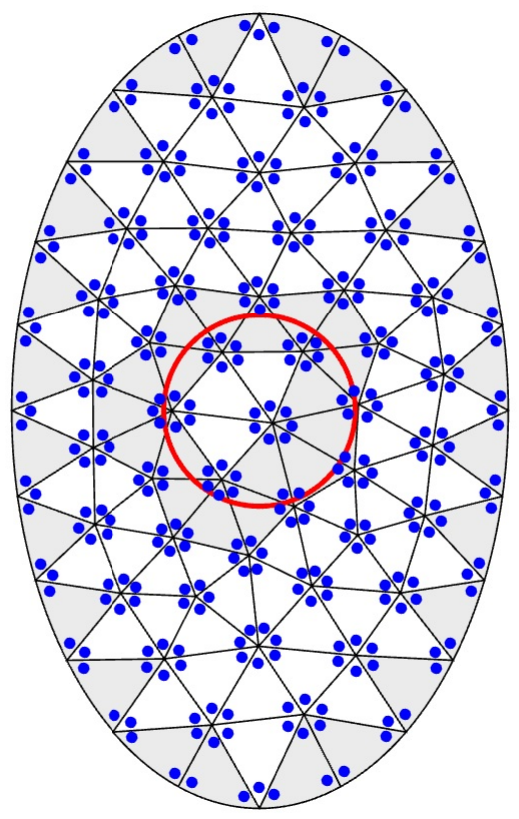}
        \caption{DG method}
        \label{fig:dof_dg}
    \end{subfigure}
    \hfill
    \begin{subfigure}{0.3\textwidth}
        \centering
        \includegraphics[width=\linewidth]{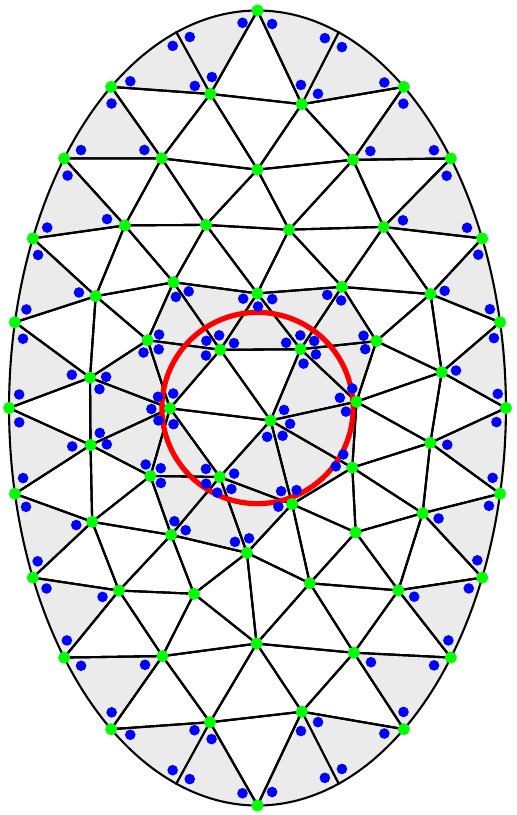}
        \caption{SDG method}
        \label{fig:dof_sdg}
    \end{subfigure}
\hfill
    % ---------- Second row ----------
    \begin{subfigure}{0.3\textwidth}
        \centering
        \includegraphics[width=\linewidth]{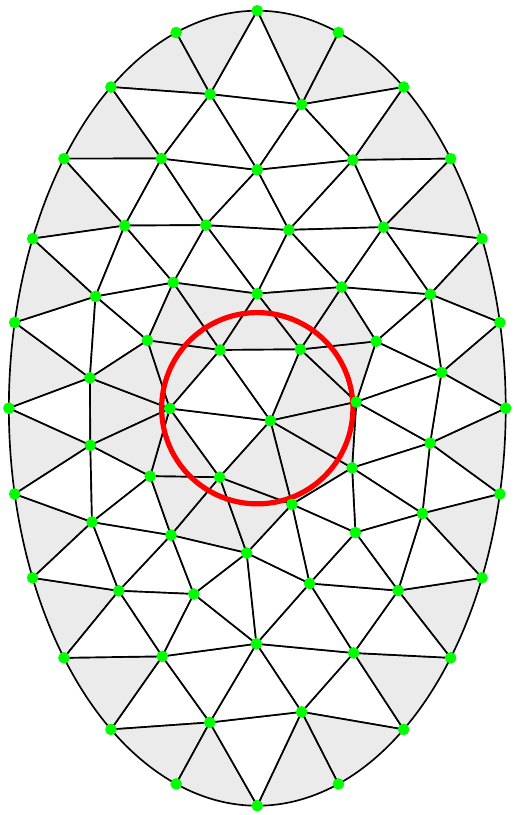}
        \caption{CG method}
        \label{fig:dof_cg}
    \end{subfigure}

    \caption{Distribution of global and local DoFs for the DG, SDG, and CG methods for $m=1$.}
    \label{fig:DOF}
\end{figure}

Compared with the SDG method of \cite{2026LiuMeghaichiZhang}, which is formulated on Cartesian meshes of polygonal domains and penalizes only the edges of interface elements, the present method treats boundary and interface elements alike: both carry non-polynomial local spaces whose traces do not match those of the neighboring Lagrange spaces, and both are therefore coupled weakly through the penalty terms on $\Es$. All remaining edges retain the CG connectivity. Figure~\ref{fig:DOF} compares the distributions of global degrees of freedom for the symmetric SDG method \eqref{SDG}, the fully DG method (each element carries its own unknowns), and the CG method (all nodes are shared) for $m=1$. A quantitative comparison is given in Section~\ref{DOF_R}.

%--------------------------------------------------------------------------
\subsection{Well-Posedness}\label{sec:wellposed}
%--------------------------------------------------------------------------
{For $v\in V_h$, define
\begin{equation}\label{4.4}
\triplenorm{v}^2_{h}:=\sum_{K\in\mathcal{T}_h}\|\sqrt\beta\nabla v\|^2_{L^2(K)}+\sum_{B\in\Es}\frac{\sigma_0}{h_B}\|\jump{v}\|^2_{L^2(B)}+\sum_{B\in\Es}\frac{h_B}{\sigma_0}\|\avg{\beta\nabla v\cdot\mathbf{n}_B}\|^2_{L^2(B)} .
\end{equation}
This is a norm on $V_h$: if $\triplenorm{v}_h=0$, then $v$ is constant on each element, continuous across every interior edge (across $\Eo\setminus\Es$ by definition of $V_h$, across $\Es\cap\Eo$ because the jumps vanish), and zero on $\Gamma_b$ because $\Eh^b\subset\Es$; hence $v=0$. 

We first prove the coercivity of $a_h(\cdot, \cdot)$ with respect to $\triplenorm{\cdot}_h$.

\begin{lemma}\label{lem:coer}
{There is a constant $C_t>0$, depending only on the constant in \eqref{3.10} and on the shape regularity, such that if $\sigma\ge4C_t$ in \eqref{eq:sigma0}, then}
\begin{equation}
{\frac{2}{5}}\triplenorm{v_h}_h^2
\leq
a_h(v_h,v_h),
\qquad
\forall v_h\in S_h^m(\Omega).
\label{3.9}
\end{equation}
\end{lemma}
\begin{proof}
{Fix $v_h\in S_h^m(\Omega)$ and abbreviate
\[
\mathsf A:=\sum_{K\in\Th}\|\sqrt\beta\nabla v_h\|^2_{L^2(K)},\qquad
\mathsf J:=\sum_{B\in\Es}\frac{\sigma_0}{h_B}\|\jump{v_h}\|^2_{L^2(B)},\qquad
\mathsf F:=\sum_{B\in\Es}\frac{h_B}{\sigma_0}\|\avg{\beta\nabla v_h\cdot\mathbf n_B}\|^2_{L^2(B)},
\]
so that $\triplenorm{v_h}_h^2=\mathsf A+\mathsf J+\mathsf F$ and 
\[a_h(v_h,v_h)=\mathsf A+\mathsf J-2\sum_{B\in\Es}\int_B\avg{\beta\nabla v_h\cdot\mathbf n_B}\jump{v_h}\,ds.\]

Since $\avg{w}^2\le\frac12\big((w|_{K_{B,1}})^2+(w|_{K_{B,2}})^2\big)$ on interior edges, each element has three edges, and $h_B\simeq h$, the trace inequality \eqref{3.10} gives
\begin{equation*}
\sum_{B\in\Es}h_B\|\avg{\beta\nabla v_h\cdot\mathbf n_B}\|^2_{L^2(B)}
\lesssim h\sum_{K\in\Th}\|\beta\nabla v_h\|^2_{L^2(\partial K)}
\le C_t\frac{\beta_{\max}^2}{\beta_{\min}}\,\mathsf A.
\end{equation*}
Hence,
\begin{equation}\label{eq4.4.9}
\mathsf F\le\frac{C_t\beta_{\max}^2}{\sigma_0\beta_{\min}}\mathsf A=\frac{C_t}{\sigma}\mathsf A .
\end{equation}

By the Cauchy--Schwarz inequality and Young's inequality $2\sqrt{\mathsf F\mathsf J}\le 2\mathsf F+\frac12\mathsf J$,
\begin{equation*}
\begin{aligned}
2\Big|\sum_{B\in\Es}\int_B\avg{\beta\nabla v_h\cdot\mathbf{n}_B}\jump{v_h}\,ds\Big|
&\le 2\sum_{B\in\Es}\Big(\frac{h_B}{\sigma_0}\Big)^{\frac12}\|\avg{\beta\nabla v_h\cdot\mathbf{n}_B}\|_{L^2(B)}\Big(\frac{\sigma_0}{h_B}\Big)^{\frac12}\|\jump{v_h}\|_{L^2(B)}\\
&\le 2\sqrt{\mathsf F\mathsf J}\le\frac{2C_t}{\sigma}\mathsf A+\frac12\mathsf J .
\end{aligned}
\end{equation*}
Therefore, for $\sigma\ge4C_t$,
\[
a_h(v_h,v_h)\ge\Big(1-\frac{2C_t}{\sigma}\Big)\mathsf A+\frac12\mathsf J\ge\frac12(\mathsf A+\mathsf J)
\]
and
\[
\triplenorm{v_h}_h^2=\mathsf A+\mathsf J+\mathsf F\le\Big(1+\frac{C_t}{\sigma}\Big)\mathsf A+\mathsf J\le\frac54(\mathsf A+\mathsf J)\le\frac52\,a_h(v_h,v_h).\qedhere
\]}
\end{proof}

The continuity of $a_h(\cdot, \cdot)$ is established in the following lemma. 

\begin{lemma}\label{lem:cont}
For all $w,v\in V_h$,
\begin{equation}
\left|a_h(w,v)\right|
\le 4
\triplenorm{w}_h\triplenorm{v}_h .
\label{4.7}
\end{equation}
\end{lemma}
\begin{proof}
{By the Cauchy--Schwarz inequality, each of the four sums in \eqref{4.2} is bounded by $\triplenorm{w}_h\triplenorm{v}_h$. For the first and last sums, this follows immediately from the definition \eqref{4.4}. For the first consistency sum, we have, as in the proof of Lemma~\ref{lem:coer},
\[
\Big|\sum_{B\in\Es}\int_B\avg{\beta\nabla w\cdot\mathbf n_B}\jump{v}\,ds\Big|
\le\Big(\sum_{B\in\Es}\frac{h_B}{\sigma_0}\|\avg{\beta\nabla w\cdot\mathbf n_B}\|^2_{L^2(B)}\Big)^{\frac12}\Big(\sum_{B\in\Es}\frac{\sigma_0}{h_B}\|\jump{v}\|^2_{L^2(B)}\Big)^{\frac12}\le\triplenorm{w}_h\triplenorm{v}_h.
\]
The estimate for the second consistency sum can be derived similarly. 
}
\end{proof}

\begin{theorem}\label{thm4.1}
Let $\sigma\ge4C_t$. Then the SDG method \eqref{SDG} has a unique solution $u_h\in S_h^m(\Omega)$.
\end{theorem}
\begin{proof}
By Lemmas~\ref{lem:coer} and \ref{lem:cont}, $a_h$ is coercive and bounded on the finite-dimensional space $S_h^m(\Omega)$ equipped with $\triplenorm{\cdot}_h$; the Lax--Milgram theorem applies.
\end{proof}

%------------------------------------------------------------------------------------------------------------------------
\subsection{A Priori Error Estimates}\label{sec:error}
%-------------------------------------------------------------------------------------------------------------------------
We first estimate the projection error $u-\Pi^m_hu$ in the energy norm. {Let $\Th^{s+}$ denote the set of elements having at least one edge in $\Es$; it consists of the selected elements and their neighbors, and $\Th^{s}\subset\Th^{s+}$.

\begin{lemma}\label{lemma4.4}
Let $u\in\mathcal{H}^{m+1}(\Omega)$ satisfy \eqref{1.4} and let $\Pi^m_h u$ be defined by \eqref{Proj}. Then
\begin{equation}
\triplenorm{u-\Pi^m_h u} _{h}\lesssim\frac{\beta_{\max}}{\sqrt{\beta_{\min}}}h^{m}\|u\|_{\mathcal{H}^{m+1}(\Omega)}.\label{4.5}
\end{equation}
\end{lemma}
\begin{proof}
Set $w:=u-\Pi^m_hu$. By Theorem~\ref{thm:approx},
\begin{equation}
\sum_{K\in\mathcal{T}_h}\|\sqrt{\beta}\nabla w\|^2_{L^2(K)}\leq \beta_{\max}\sum_{K\in\mathcal{T}_h}\|\nabla w\|^2_{L^2(K)}\lesssim\beta_{\max} h^{2m}\|u\|^2_{\mathcal{H}^{m+1}(\Omega)}.\label{4.2.10}
\end{equation}
{Since $w|_K\in\Hc^2(K)\subset H^1(K)$ for every $K$, the trace inequality $\|w\|^2_{L^2(\partial K)}\lesssim h^{-1}\|w\|^2_{L^2(K)}+h\|\nabla w\|^2_{L^2(K)}$ (valid on straight and on curved boundary elements) and Theorem~\ref{thm:approx} with $j=0,1$ give}
\begin{equation}\label{4.2.11}
\begin{aligned}
\sum_{B\in\Es}\frac{\sigma_0}{h_B}\|\jump{w}\|^2_{L^2(B)}
&\lesssim \frac{\sigma_0}{h}\sum_{K\in\Th^{s+}}\| w\|^2_{L^2(\partial K)}\\
&\lesssim \frac{\sigma_0}{h}\sum_{K\in\Th^{s+}}\Big(\frac{\|w\|^2_{L^2(K)}}{h} +h\|\nabla w\|^2_{L^2(K)}  \Big)
\lesssim \sigma_0 h^{2m}\|u\|^2_{\mathcal{H}^{m+1}(\Omega)}.
\end{aligned}
\end{equation}
{For the flux averages we use Theorem~\ref{thm:traceH2} on interface elements and the standard trace inequality for $\beta\nabla w\in H^1(K)^2$ on the other elements of $\Th^{s+}$ (where $\beta$ is constant); in both cases
\[
\|\beta\nabla w\|^2_{L^2(\partial K)}\lesssim\Big(\frac{\beta_{\max}}{\beta_{\min}}\Big)^2\Big(h^{-1}\|\beta\nabla w\|^2_{L^2(K)}+h\|\beta D^2w\|^2_{L^2(K^-\cup K^+)}\Big).
\]
Hence, by Theorem~\ref{thm:approx} with $j=1,2$,}
\begin{equation}
\begin{aligned}
\sum_{B\in\Es}\frac{h_B}{\sigma_0}\|\avg{\beta\nabla w\cdot \mathbf{n}_B}\|^2_{L^2(B)}&\lesssim \frac{h}{\sigma_0}\sum_{K\in\Th^{s+}}\|\beta\nabla w\|^2_{L^2(\partial K)}\\
&\lesssim \frac{h}{\sigma_0}\Big(\frac{\beta_{\max}}{\beta_{\min}}\Big)^2\sum_{K\in\Th^{s+}}\Big(h^{-1}\|\beta\nabla w\|^2_{L^2(K)} +h\|\beta D^2 w\|^2_{L^2(K^-\cup K^+)}  \Big)\\
&\lesssim \frac{\beta_{\max}^4}{\sigma_0\beta_{\min}^2} h^{2m}\|u\|^2_{\mathcal{H}^{m+1}(\Omega)}.\label{4.2.11.2}
\end{aligned}
\end{equation}
Combining \eqref{4.2.10}--\eqref{4.2.11.2} and recalling $\sigma_0=\sigma\beta_{\max}^2/\beta_{\min}$ with $\sigma$ a fixed constant,
\begin{equation*}
\triplenorm{w}^2_{h}\lesssim\left(\beta_{\max}+\sigma_0+\frac{\beta_{\max}^4}{\sigma_0\beta_{\min}^2}\right)h^{2m}\|u\|^2_{\mathcal{H}^{m+1}(\Omega)}\lesssim\frac{\beta_{\max}^2}{\beta_{\min}}h^{2m}\|u\|^2_{\mathcal{H}^{m+1}(\Omega)},
\end{equation*}
where we used $\beta_{\max}\le\beta_{\max}^2/\beta_{\min}$.
\end{proof}

\begin{theorem}\label{thm:energy}
Let $u\in\mathcal H^{m+1}(\Omega)$ be the solution of \eqref{1.1}--\eqref{1.4} and let $u_h\in S_h^m(\Omega)$ be the SDG solution determined by \eqref{SDG} with $\sigma\ge4C_t$. Then
\begin{equation}\label{4.2.17}
\triplenorm{u-u_h}_h
\lesssim
\frac{\beta_{\max}}{\sqrt{\beta_{\min}}}
h^m
\|u\|_{\mathcal H^{m+1}(\Omega)}.
\end{equation}
\end{theorem}
\begin{proof}
Let $v_h=u_h-\Pi^m_hu\in S_h^m(\Omega)$. By Lemma~\ref{lem:coer}, the Galerkin orthogonality \eqref{eq:GO}, and Lemma~\ref{lem:cont},
\begin{equation*}
    \frac{2}{5}\triplenorm{v_h}_h^2\leq a_h(v_h,v_h)= a_h(u-\Pi^m_h u,v_h)\leq 4\triplenorm{u-\Pi^m_h u}_h\triplenorm{v_h}_h ,
\end{equation*}
so $\triplenorm{u_h-\Pi^m_h u}_h\leq 10\triplenorm{u-\Pi^m_h u}_h$. The triangle inequality and Lemma~\ref{lemma4.4} give
\[
\triplenorm{u-u_h}_h\leq \triplenorm{u-\Pi^m_hu}_h+\triplenorm{u_h-\Pi^m_hu}_h\leq 11\triplenorm{u-\Pi^m_hu}_h\lesssim\frac{\beta_{\max}}{\sqrt{\beta_{\min}}}h^m\|u\|_{\mathcal{H}^{m+1}(\Omega)}.\qedhere
\]
\end{proof}

For the $L^2$ estimate, we use a duality argument. For $\theta\in L^2(\Omega)$, let $w=w_\theta$ be the solution of the auxiliary problem \eqref{1.1}--\eqref{1.3} with $f=\theta$ and $g=0$. The argument requires regularity of $w$ and an approximation estimate for $w$ in the energy norm. The latter does not follow from Lemma~\ref{lemma4.4}, which is stated for functions in $\Hc^{m+1}(\Omega)$ satisfying the extended conditions \eqref{1.4}, whereas $w$ is in general only in $\Hc^2(\Omega)$. The following assumption states the required dual regularity.

\begin{assumption}\label{ass:dual}
For every $\theta\in L^2(\Omega)$, the solution $w$ of the auxiliary problem satisfies $w\in\Hc^2(\Omega)$ with
\begin{equation}\label{eq:dualreg}
\beta^-\|w\|_{H^2(\Omega^-)}+\beta^+\|w\|_{H^2(\Omega^+)}\lesssim\|\theta\|_{L^2(\Omega)}.
\end{equation}
%and the projection $\Pi_hw\in S_h^m(\Omega)$ defined by \eqref{Proj} satisfies
%\begin{equation}\label{eq:dualapprox}
%\triplenorm{w-\Pi_hw}_h\lesssim\frac{\beta_{\max}}{\sqrt{\beta_{\min}}}\,h\,\|w\|_{\Hc^2(\Omega)} .
%\end{equation}
\end{assumption}

The regularity estimate \eqref{eq:dualreg} is the classical one for interface problems with smooth interfaces; see \cite{1998ChenZou} and the references therein. %For $m=1$, the estimate \eqref{eq:dualapprox} is exactly Lemma~\ref{lemma4.4}, since \eqref{1.4} is void in that case. For $m\ge2$, \eqref{eq:dualapprox} is a low-regularity approximation property of the degree-$m$ GC-IFE space, which can be established by repeating the arguments of Section~\ref{sec:approx} with $m+1$ replaced by $2$ in Lemmas~\ref{lemma3.1}--\ref{lemma3.3}, as done for rectangular Frenet elements in \cite{2024AdjeridLinMeghaichi,2025AdjeridLinMeghaichi}; we do not pursue this here

\begin{theorem}\label{thm:L2}
Under the assumptions of Theorem~\ref{thm:energy} and Assumption~\ref{ass:dual},
\begin{equation}\label{eq: L2}
\|u-u_h\|_{L^2(\Omega)}\lesssim\left(\frac{\beta_{\max}}{\beta_{\min}}\right)^2h^{m+1}\|u\|_{\mathcal{H}^{m+1}(\Omega)}.
\end{equation}
\end{theorem}
%\begin{proof}
%Let $\theta:=u-u_h\in V_h$ and let $w$ be the corresponding solution of the auxiliary problem. {Since $w\in\Hc^2(\Omega)$ and $g=0$, the consistency identity \eqref{4.1} applied to $w$ gives $a_h(w,\theta)=(\theta,\theta)_{L^2(\Omega)}$}, and the Galerkin orthogonality \eqref{eq:GO} together with the symmetry of $a_h$ gives $a_h(\Pi_hw,\theta)=0$. Hence, by Lemma~\ref{lem:cont},
%\begin{equation*}
%\|\theta\|_{L^2(\Omega)}^2=a_h(w-\Pi_h w,\theta)\leq  4\triplenorm{w-\Pi_hw}_{h}\triplenorm{\theta}_{h}.
%\end{equation*}
%By \eqref{eq:dualapprox}, \eqref{eq:dualreg} and $\|w\|_{\Hc^2(\Omega)}\le\beta_{\min}^{-1}(\beta^-\|w\|_{H^2(\Omega^-)}+\beta^+\|w\|_{H^2(\Omega^+)})$,
%\begin{equation}
%\triplenorm{w-\Pi_h w}_h\lesssim \frac{\beta_{\max}}{\sqrt{\beta_{\min}}}h\|w\|_{\mathcal{H}^{2}(\Omega)}\lesssim \frac{\beta_{\max}}{\beta_{\min}^{3/2}}h\|\theta\|_{L^2(\Omega)}.\label{4.19}
%\end{equation}
%Combining this with Theorem~\ref{thm:energy} yields $\|\theta\|_{L^2(\Omega)}^2\lesssim (\beta_{\max}/\beta_{\min})^2h^{m+1}\|u\|_{\mathcal{H}^{m+1}(\Omega)}\|\theta\|_{L^2(\Omega)}$, which is \eqref{eq: L2}.
%\end{proof}
\begin{proof}
Let $\theta:=u-u_h\in V_h$ and let $w$ be the corresponding solution of the auxiliary problem. {Since $w\in\Hc^2(\Omega)$ and $g=0$, the consistency identity \eqref{4.1} applied to $w$ gives $a_h(w,\theta)=(\theta,\theta)_{L^2(\Omega)}$}, and the Galerkin orthogonality \eqref{eq:GO} together with the symmetry of $a_h$ gives $a_h(\Pi_h^1w,\theta)=0$. Hence, by Lemma~\ref{lem:cont} and Theorem~\ref{thm:energy}, we have

\begin{equation}\label{thetaL2}
\|\theta\|_{L^2(\Omega)}^2=a_h(w-\Pi^1_h w,\theta)\leq  4\triplenorm{w-\Pi^1_hw}_{h}\triplenorm{\theta}_{h}\lesssim \frac{\beta_{\max}}{\sqrt{\beta_{\min}}}
h^m
\|u\|_{\mathcal H^{m+1}(\Omega)}\triplenorm{w-\Pi^1_hw}_{h}.
\end{equation}
By \eqref{eq:dualreg} and $\|w\|_{\Hc^2(\Omega)}\le\beta_{\min}^{-1}(\beta^-\|w\|_{H^2(\Omega^-)}+\beta^+\|w\|_{H^2(\Omega^+)})$, we get

\begin{equation}
\triplenorm{w-\Pi^1_h w}_h\lesssim \frac{\beta_{\max}}{\sqrt{\beta_{\min}}}h\|w\|_{\mathcal{H}^{2}(\Omega)}\lesssim \frac{\beta_{\max}}{\beta_{\min}^{3/2}}h\|\theta\|_{L^2(\Omega)}.\label{4.19}
\end{equation}
Combining \eqref{thetaL2} and \eqref{4.19} leads to $\|\theta\|_{L^2(\Omega)}^2\lesssim (\frac{\beta_{\max}}{\beta_{\min}})^2h^{m+1}\|u\|_{\mathcal{H}^{m+1}(\Omega)}\|\theta\|_{L^2(\Omega)}$. Dividing by $\|\theta\|_{L^2(\Omega)}$ gives \eqref{eq: L2}.
\end{proof}

\begin{remark}\label{rem:special}
 If $\Omega$ is a polygon and $\Th$ is exactly fitted to $\partial\Omega$, then $\Th^b=\varnothing$, $\Th^s=\Th^i$, and all results above hold verbatim; the method then reduces to the SDG method of \cite{2026LiuMeghaichiZhang} with the triangular GC-IFE spaces of \cite{2025LinLinZhang}, for which the present analysis provides the error estimates.
 \end{remark}
 \begin{remark}\label{rem:special2}
  Replacing $\Es$ by $\Eh$ and $S_h^m(\Omega)$ by $\widetilde S_h^m(\Omega)$ in \eqref{4.2}--\eqref{SDG} gives the symmetric interior penalty DG method on the broken space $\widetilde{\Hc}^{2}(\Th)$. Lemmas~\ref{lem:coer}--\ref{lem:cont} and Theorems~\ref{thm:energy}--\ref{thm:L2} remain valid with the same proofs, the only change being that the sums in \eqref{4.2.11}--\eqref{4.2.11.2} extend over all elements.
\end{remark}

\section{Numerical Experiments} \label{sec5}
In this section, we present numerical examples to demonstrate the performance of the GC-HIFE space and the SDG method. The first example is posed on a square domain with a circular interface, so that only interface elements are in $\mathcal T_h^{s}$. It is used to quantify the reduction in degrees of freedom achieved by the SDG formulation and to verify the convergence rates established in Theorems~\ref{thm:energy}--\ref{thm:L2} for large coefficient contrasts. The second example is posed on an elliptical domain with a flower-shaped interface, so that both boundary and interface elements are selected to be treated discontinuously. In all tests, the penalty parameter is $\sigma_0=\beta_{\max}^2/\beta_{\min}$, i.e., $\sigma=1$ in \eqref{eq:sigma0}; no instability was observed for this choice.

\subsection{Square Domain with a Circular Interface}\label{sec:rect}
We take $\Omega=(-1,1)^2$ and the level-set function $\phi(x,y)=x^2+y^2-r^2$ {with $r=0.44$}. The circular interface $\Gamma_i=\{\phi=0\}$ separates $\Omega^-=\{\phi<0\}$ from $\Omega^+=\{\phi>0\}\cap\Omega$. The exact solution is
\begin{equation}\label{eq:5.1}
u(x,y)=
\begin{cases}
\displaystyle \frac{\cos(\phi(x,y))}{\beta^-},
& (x,y)\in\Omega^-,\\[6pt]
\displaystyle \frac{\cos(\phi(x,y))}{\beta^+}
+\frac{1}{\beta^-}-\frac{1}{\beta^+},
& (x,y)\in\Omega^+.
\end{cases}
\end{equation}
Figure~\ref{fig:mesh_rectangular} shows an example of the interface-unfitted triangular mesh. We use uniform triangular meshes obtained by dividing $\Omega$ into $N\times N$ squares, each cut into two triangles, so that $h=O(N^{-1})$.

\begin{figure}[H]
    \centering
    \begin{subfigure}[b]{0.5\linewidth}
        \centering\includegraphics[height=5cm]{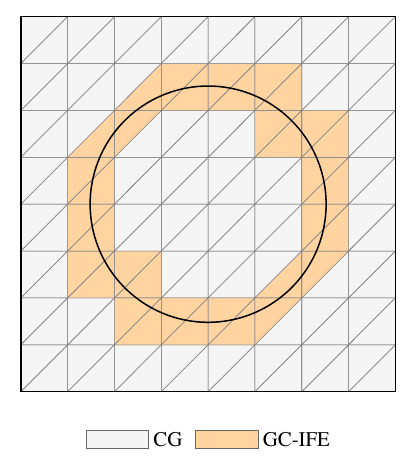}
         \caption{Square domain.}\label{fig:mesh_rectangular}
    \end{subfigure}\hfill
    \begin{subfigure}[b]{0.5\linewidth}
        \centering\includegraphics[height=5cm]{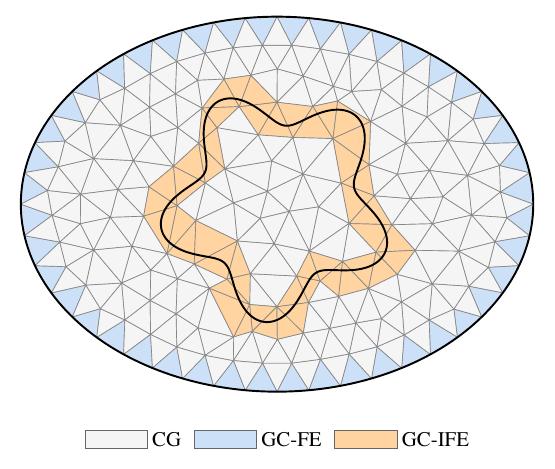}
         \caption{Elliptical domain.}\label{fig:elliptical-mesh}
    \end{subfigure}
    \caption{Meshes and element types: (a) square domain with a circular interface; (b) elliptical domain with a flower-shaped interface.
    }
    \label{fig:mesh-examples}
\end{figure}

\subsubsection{Degrees of Freedom}\label{DOF_R}
Let $\mathrm{DoF}_D$, $\mathrm{DoF}_S$, and $\mathrm{DoF}_C$ denote the numbers of global degrees of freedom of the DG, SDG, and CG methods of degree $m$. {On the meshes described above, $\mathrm{DoF}_D(m)=d_m\cdot 2N^2$ and $\mathrm{DoF}_C(m)=(mN+1)^2$. In the SDG method, each interface element carries private copies of its $3m$ boundary nodes, while the $m-1$ interior nodes of each cut edge $B\in\Eo^{i}$ (the interior edges intersected by $\Gamma_i$) are no longer needed as shared nodes; for the meshes used here, where $\Gamma_i$ crosses every interface element exactly once, $|\Eo^{i}|=|\Th^i|$ and
\begin{equation*}
\mathrm{DoF}_S(m)=\mathrm{DoF}_C(m)+3m|\Th^i|-(m-1)|\Eo^{i}|=(mN+1)^2+(2m+1)|\Th^i| .
\end{equation*}
Since $|\Th^i|=O(N)$, we obtain, as $N\to\infty$,
\begin{equation*}
\begin{aligned}
&\frac{\mathrm{DoF}_{D}(m)}{\mathrm{DoF}_{S}(m)}=\frac{(m+1)(m+2)N^2}{(mN+1)^2+(2m+1)|\mathcal{T}_h^i|}\to \frac{(m+1)(m+2)}{m^2},
\\
&\frac{\mathrm{DoF}_{C}(m)}{\mathrm{DoF}_{S}(m)}=\frac{(mN+1)^2}{(mN+1)^2+(2m+1)|\mathcal{T}_h^i|}\to 1 .
\end{aligned}
\end{equation*}}

Thus, $\mathrm{DoF}_S$ is asymptotically smaller than $\mathrm{DoF}_D$, while remaining asymptotically equivalent to $\mathrm{DoF}_C$.   Figure~\ref{fig:DoF_asymtotic} and Table~\ref{tab3} illustrate that, on the same mesh, the proposed SDG method uses fewer degrees of freedom than the full DG method, with counts close to those of the CG method when the mesh is refined.

Moreover, comparing the numbers of global degrees of freedom of the DG method of degree $m$ and the SDG
method of degree $m+1$, we obtain
\begin{equation}
\frac{\mathrm{DoF}_{D}(m)}{\mathrm{DoF}_{S}(m+1)}
=
\frac{(m+1)(m+2)N^2}
{\bigl((m+1)N+1\bigr)^2+(2m+3)|\mathcal{T}_h^i|}
\longrightarrow
\frac{m+2}{m+1},
\qquad \text{as } N\to\infty.
\end{equation}
Thus, for fixed $m$ and sufficiently fine meshes, the SDG method
of degree $m+1$ uses fewer global degrees of freedom than the
DG method of degree $m$. Figure~\ref{fig:DoF_comparison} shows that the DG method of degree $m$ requires more global degrees of freedom than the SDG method of degree $m+1$ for $1\le m\le 5$.

\commentout{
        \begin{table}[h]
        \centering
        \begin{tabular}{cccccccccc c}
        \toprule
        & & & \multicolumn{2}{c}{$m=1$} & \multicolumn{2}{c}{$m=2$} & \multicolumn{2}{c}{$m=3$}& \multicolumn{2}{c}{$m=4$}\\
        \cmidrule(lr){4-5} \cmidrule(lr){6-7}\cmidrule(lr){8-9}\cmidrule(lr){10-11}
        $N$ & $|\mathcal{T}_h^i|$ & $|\mathring{\mathcal{E}}_h^i|$
        & $\frac{\mathrm{DoF}_{D}}{\mathrm{DoF}_{S}}$  
        & $\frac{\mathrm{DoF}_C}{\mathrm{DoF}_{S}}$ 
        & $\frac{\mathrm{DoF}_{D}}{\mathrm{DoF}_{S}}$  
        & $\frac{\mathrm{DoF}_C}{\mathrm{DoF}_{S}}$
        & $\frac{\mathrm{DoF}_{D}}{\mathrm{DoF}_{S}}$  
        & $\frac{\mathrm{DoF}_C}{\mathrm{DoF}_{S}}$ 
        & $\frac{\mathrm{DoF}_{D}}{\mathrm{DoF}_{S}}$  
        & $\frac{\mathrm{DoF}_C}{\mathrm{DoF}_{S}}$  \\
        \midrule
        8    & 34   & 34   & 2.0984 & 0.4426 & 1.6732 & 0.6296 & 1.4832 & 0.7242 & 1.3763 & 0.7806\\
        16   & 62   & 62   & 3.2337 & 0.6084 & 2.1959 & 0.7784 & 1.8060 & 0.8469 & 1.6057 & 0.8833\\
        32   & 130  & 130  & 4.1542 & 0.7363 & 2.5206 & 0.8667 & 1.9847 & 0.9118 & 1.7248 & 0.9343\\
        64   & 254  & 254  & 4.9280 & 0.8472 & 2.7442 & 0.9291 & 2.0991 & 0.9544 & 1.7982 & 0.9665\\
        128  & 502  & 502  & 5.4171 & 0.9170 & 2.8677 & 0.9634 & 2.1595 & 0.9768 & 1.8362 & 0.9831\\
        256  & 1006 & 1006 & 5.6933 & 0.9563 & 2.9323 & 0.9812 & 2.1904 & 0.9882 & 1.8554 & 0.9915\\
        512  & 2018 & 2018 & 5.8422 & 0.9775 & 2.9657 & 0.9905 & 2.2061 & 0.9941 & 1.8651 & 0.9957\\
        1024 & 4038 & 4038 & 5.9200 & 0.9886 & 2.9827 & 0.9952 & 2.2141 & 0.9970 & 1.8700 & 0.9978\\
        \bottomrule
        \end{tabular}
        \caption{Ratios of the numbers of global degrees of freedom of the DG, SDG, and CG methods for $m=1,2,3,4$.}
        \label{tab3}
        \end{table}
}

\begin{figure}[H]
    \centering
    \begin{subfigure}[b]{0.5\linewidth}
        \centering\includegraphics[height=6.5cm]{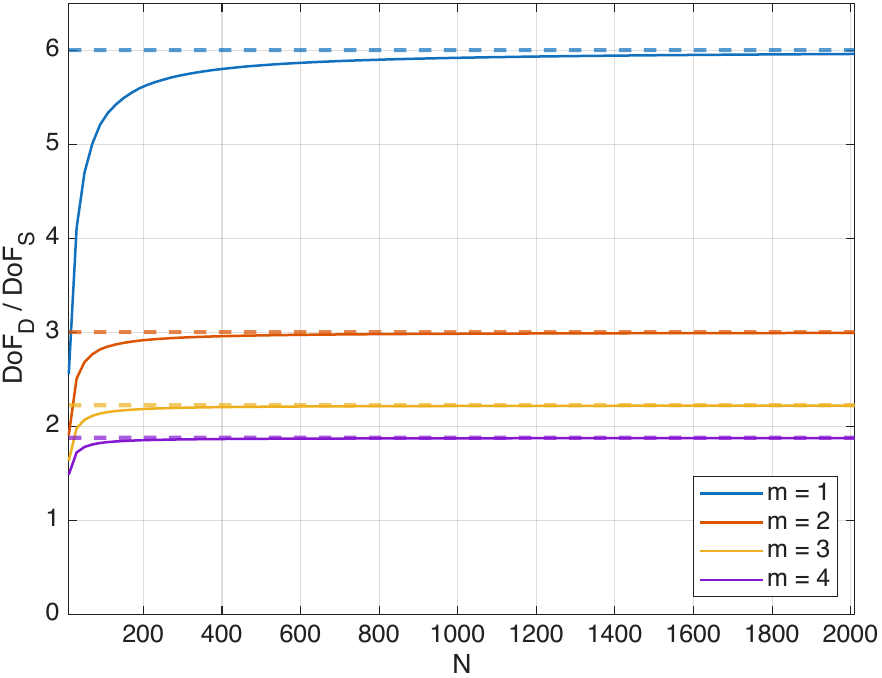}
         \caption{$\frac{\mathrm{DoF}_{D}}{\mathrm{DoF}_{S}}$}\label{fig:DG_SDG}
    \end{subfigure}\hfill
    \begin{subfigure}[b]{0.5\linewidth}
        \centering\includegraphics[height=6.5cm]{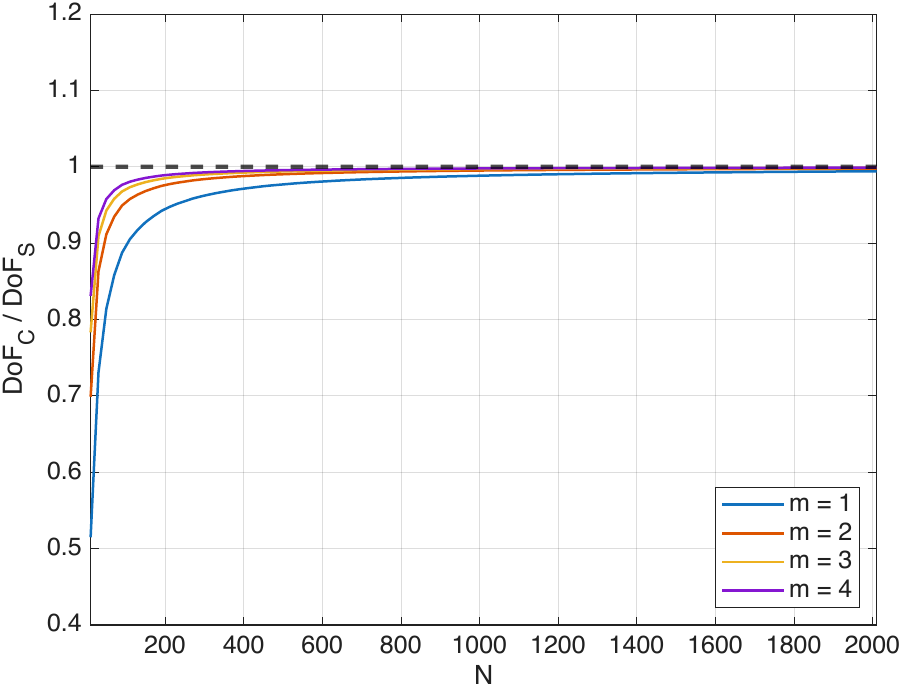}
         \caption{$\frac{\mathrm{DoF}_{C}}{\mathrm{DoF}_{S}}$}\label{fig:CG_SDG}
    \end{subfigure}
    \caption{Ratios of the degrees of freedom (DoFs)
of the DG and CG methods to those of the proposed
SDG method under mesh refinement for polynomial
degrees $1\le m\le 4$.
The dashed lines indicate the corresponding
theoretical asymptotic limits.}
    \label{fig:DoF_asymtotic}
\end{figure}
\begin{table}[ht]
\centering
\begin{tabular}{cccccccccc c}
\toprule
& & & \multicolumn{2}{c}{$m=1$} & \multicolumn{2}{c}{$m=2$} & \multicolumn{2}{c}{$m=3$}& \multicolumn{2}{c}{$m=4$}\\
\cmidrule(lr){4-5} \cmidrule(lr){6-7}\cmidrule(lr){8-9}\cmidrule(lr){10-11}
$N$ & $|\mathcal{T}_h^i|$ & $|\mathring{\mathcal{E}}_h^i|$
& $\frac{\mathrm{DoF}_{D}}{\mathrm{DoF}_{S}}$  
& $\frac{\mathrm{DoF}_C}{\mathrm{DoF}_{S}}$ 
& $\frac{\mathrm{DoF}_{D}}{\mathrm{DoF}_{S}}$  
& $\frac{\mathrm{DoF}_C}{\mathrm{DoF}_{S}}$
& $\frac{\mathrm{DoF}_{D}}{\mathrm{DoF}_{S}}$  
& $\frac{\mathrm{DoF}_C}{\mathrm{DoF}_{S}}$ 
& $\frac{\mathrm{DoF}_{D}}{\mathrm{DoF}_{S}}$  
& $\frac{\mathrm{DoF}_C}{\mathrm{DoF}_{S}}$  \\
\midrule
8    & 34   & 34   & 2.0984 & 0.4426 & 1.6732 & 0.6296 & 1.4832 & 0.7242 & 1.3763 & 0.7806\\
128  & 502  & 502  & 5.4171 & 0.9170 & 2.8677 & 0.9634 & 2.1595 & 0.9768 & 1.8362 & 0.9831\\
2048 & 8078  & 8078 & 5.9597 & 0.9943 & 2.9913 & 0.9976 & 2.2182 & 0.9985 & 1.8725 & 0.9989\\
\bottomrule
\end{tabular}
\caption{Ratios of the numbers of global degrees of freedom of the DG, SDG, and CG methods for $m=1,2,3,4$.}
\label{tab3}
\end{table}
\begin{figure}[ht]
    \centering
\includegraphics[width=0.7\textwidth]{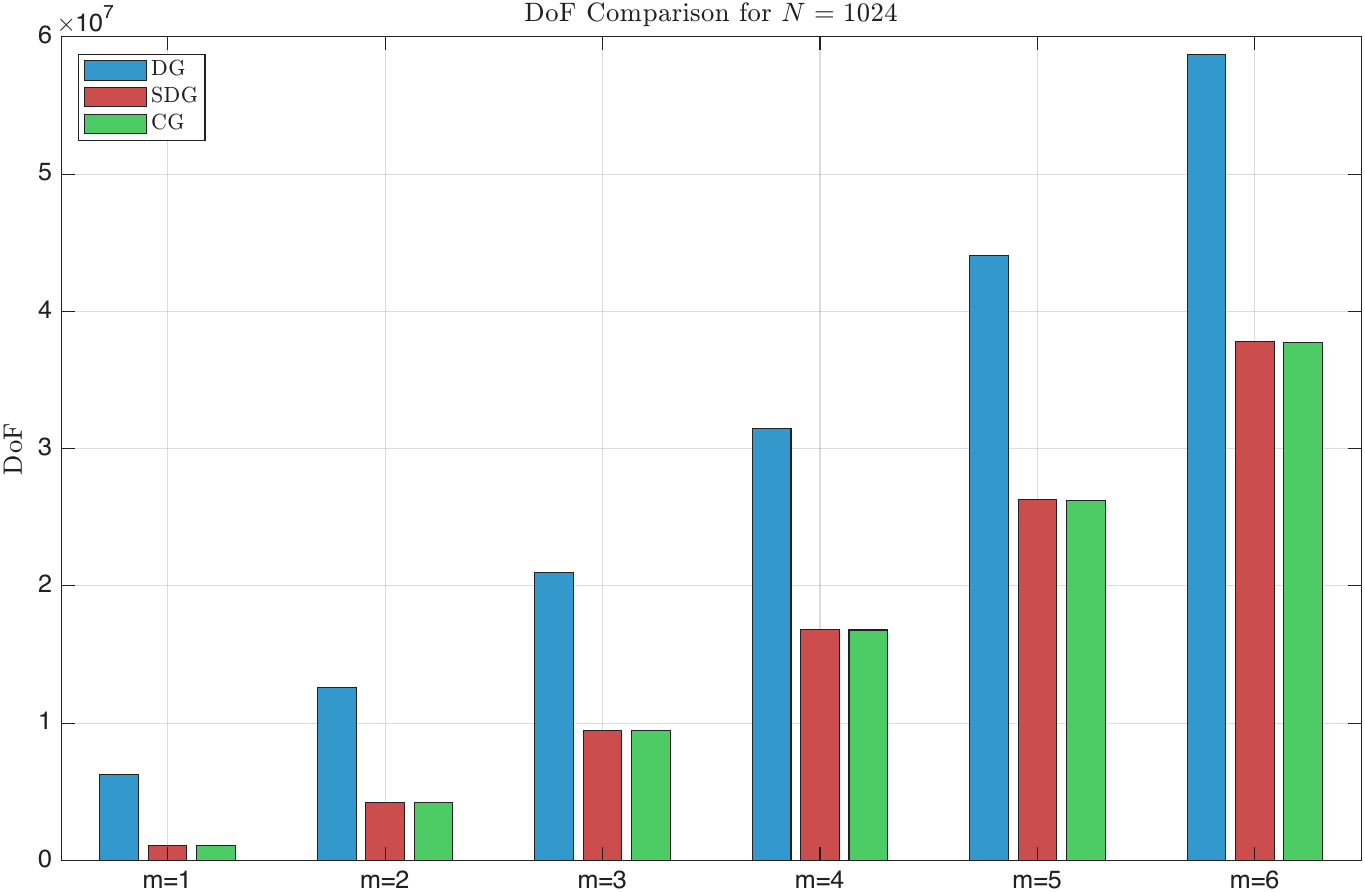} 
    \caption{Global degrees of freedom of the DG, SDG, and CG methods for $1\le m\le 6$ with $N=1024$.}
    \label{fig:DoF_comparison}
\end{figure}

\subsubsection{Convergence}
Figure~\ref{fig:R_errors} reports the $L^2$ and broken $H^1$ errors of the SDG solutions for $m=1,2,3,4$ and the coefficient pairs $(\beta^-,\beta^+)=(1,10)$ and $(1,1000)$. The errors converge at the optimal rates $O(h^{m+1})$ in $L^2$ and $O(h^m)$ in the broken $H^1$ seminorm predicted by Theorems~\ref{thm:energy}--\ref{thm:L2}, and the curves for the two contrasts are almost indistinguishable, in agreement with the explicit $\beta$-dependence of the estimates, which involves only the ratio $\beta_{\max}/\beta_{\min}$ through the penalty scaling.

\begin{figure}[htbp]
    \centering

    % First row
    \begin{subfigure}{0.48\textwidth}
        \centering
        \includegraphics[width=\linewidth]{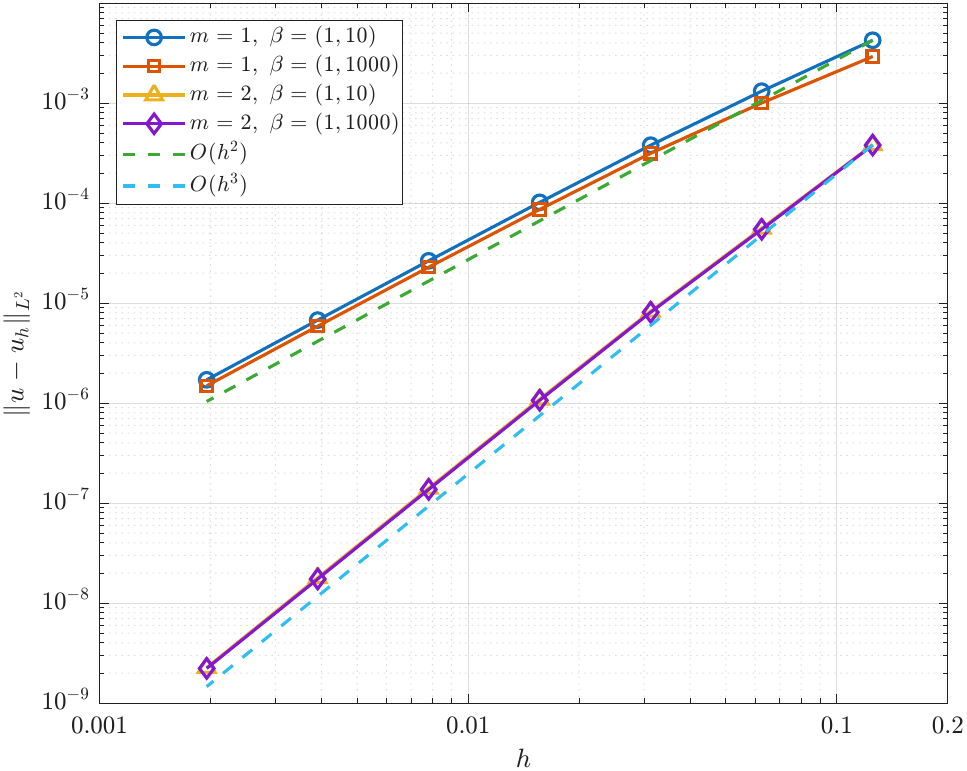}
        \caption{$L^2$ error for  $m=1,2$.}
        \label{fig:m1_L2}
    \end{subfigure}
    \hfill
    \begin{subfigure}{0.48\textwidth}
        \centering
        \includegraphics[width=\linewidth]{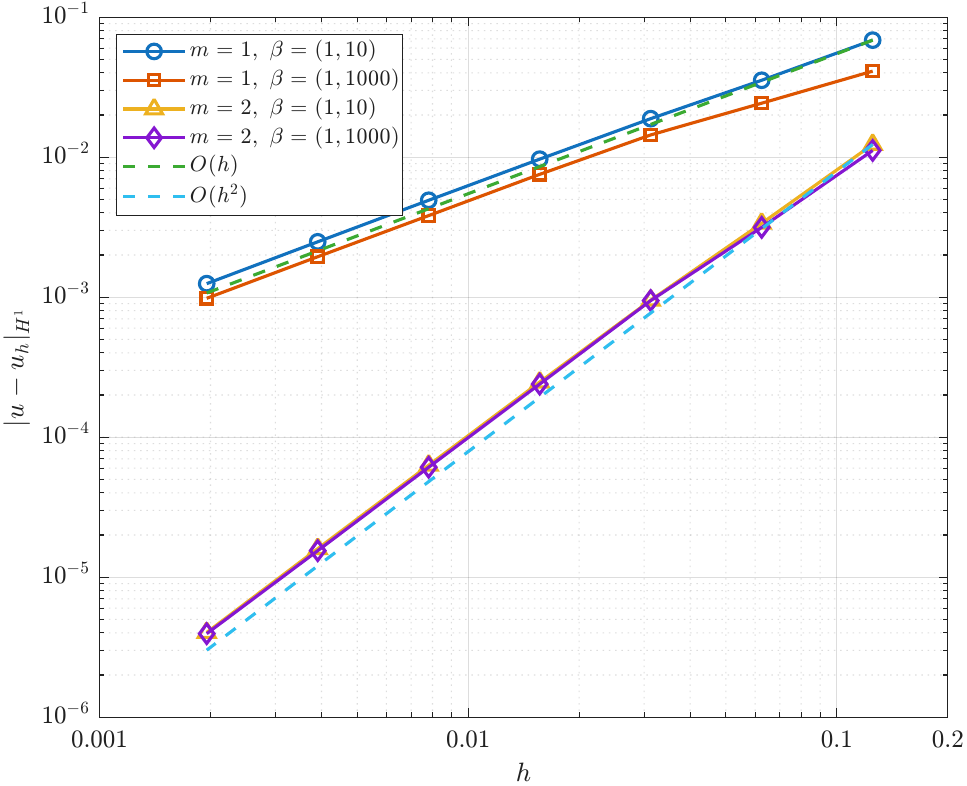}
        \caption{$H^1$-seminorm error for $m=1,2$.}
        \label{fig:m1_H1}
    \end{subfigure}

    \vspace{0.5cm}

    % Second row
    \begin{subfigure}{0.48\textwidth}
        \centering
        \includegraphics[width=\linewidth]{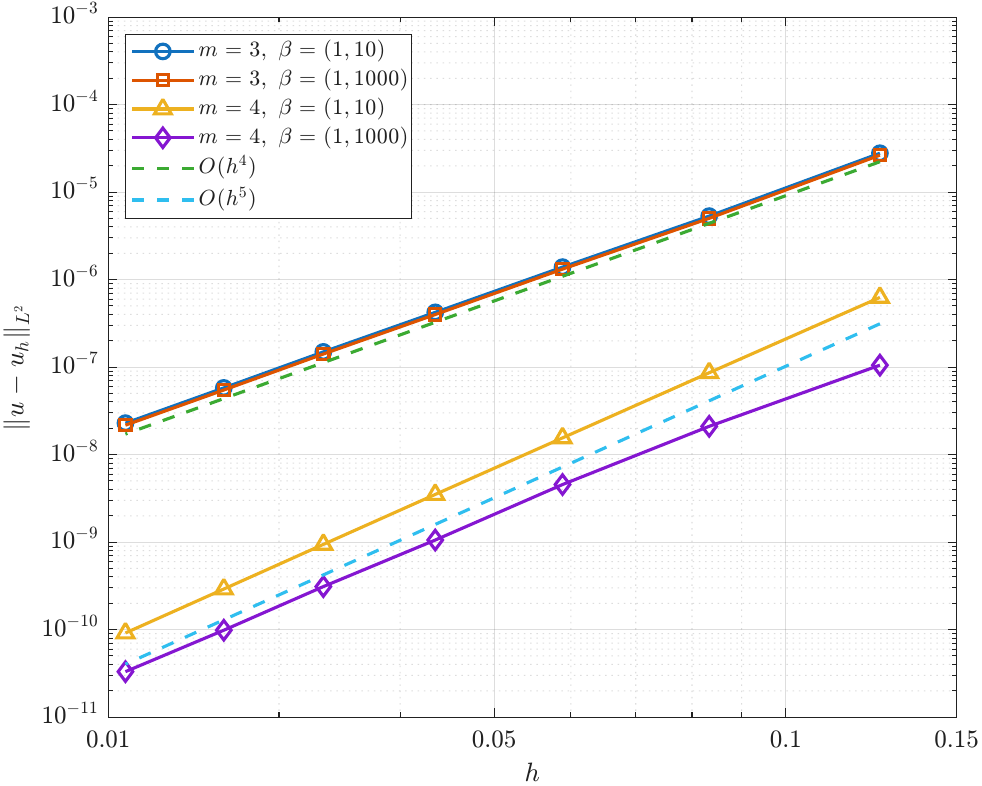}
        \caption{$L^2$ error for $m=3,4$.}
        \label{fig:m2_L2}
    \end{subfigure}
    \hfill
    \begin{subfigure}{0.48\textwidth}
        \centering
        \includegraphics[width=\linewidth]{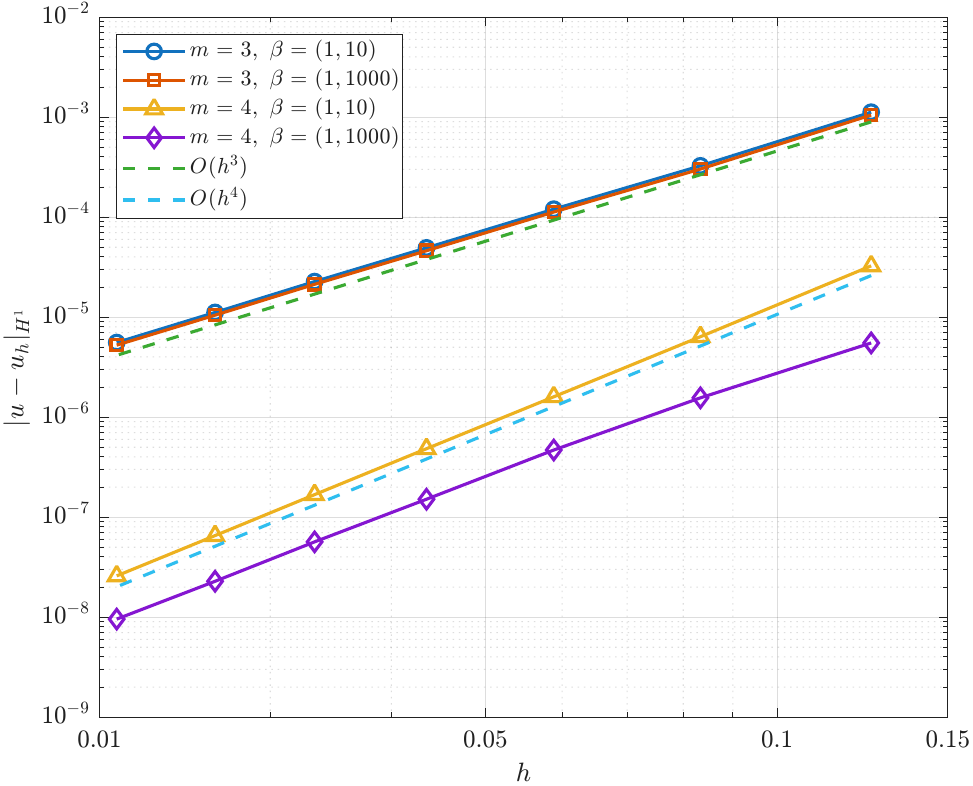}
        \caption{$H^1$-seminorm error for $m=3,4$.}
        \label{fig:m2_H1}
    \end{subfigure}

    \caption{
        Convergence results for the SDG
        method with different coefficient contrasts.
    }
    \label{fig:R_errors}
\end{figure}

\subsection{Elliptical Domain with a Flower-Shaped Interface}
\label{sec:new-flower}
We consider the elliptical domain
\begin{equation}\label{eq:new-ellipse}
 \Omega=\left\{(x,y):\frac{x^2}{2.2^2}+\frac{y^2}{1.45^2}<1\right\}.
\end{equation}
Let \(r=\sqrt{x^2+y^2}\) and \(\theta=\operatorname{atan2}(y,x)\). The interface is
\begin{equation}\label{eq:new-flower-curve}
 \Gamma_i=\{r=\rho(\theta)\},\qquad
 \rho(\theta)=1+0.2\sin\left(5\left(\theta-\frac{\pi}{6}\right)\right),
\end{equation}
described by the level-set function
\begin{equation}\label{eq:new-flower-phi}
 \phi(x,y)=\frac{1}{400}\left[\left(\frac{r}{\rho(\theta)}\right)^6-1\right]\quad(r>0),
 \qquad \phi(0,0)=-\frac{1}{400},
\end{equation}
with \(\Omega^-=\{\phi<0\}\) and \(\Omega^+=\Omega\cap\{\phi>0\}\). The exact solution is given by \eqref{eq:5.1} with this $\phi$, and \(f\) and \(g=u|_{\partial\Omega}\) are determined accordingly. The mesh is fitted to \(\partial\Omega\) and unfitted to \(\Gamma_i\), as illustrated in Figure~\ref{fig:elliptical-mesh}; {both GC-FE boundary elements and GC-IFE interface elements are present, and the Dirichlet condition is imposed weakly on the curved boundary.} Figure~\ref{fig:new-flower-errors} shows the SDG errors for \(m=1,2,3,4\) and \((\beta^-,\beta^+)=(1,10),(1,1000)\). The least-squares fitted rates are close to \(m+1\) in the \(L^2\) norm and to \(m\) in the broken \(H^1\) seminorm, {confirming that the geometry-conforming treatment of the curved boundary does not degrade the high-order accuracy of the interface approximation.}

\begin{figure}[H]
 \centering
 \includegraphics[width=0.98\linewidth]{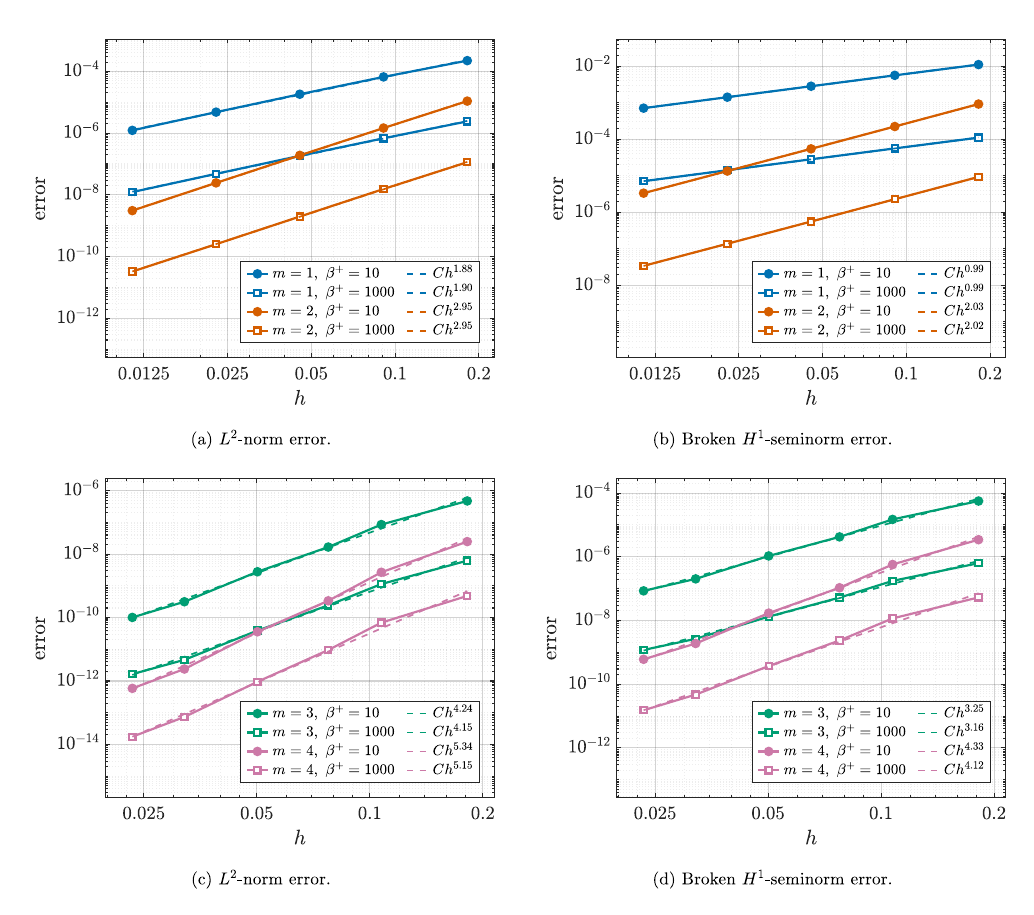}
 \caption{Elliptical domain: SDG errors for \(m=1,2,3,4\) and \((\beta^-,\beta^+)=(1,10),(1,1000)\). Dashed lines show the least-squares fits \(Ch^p\).}
 \label{fig:new-flower-errors}
\end{figure}

\section{Conclusion}\label{sec6}
We have developed and analyzed a high-order geometry-conforming selective discontinuous Galerkin method for elliptic interface problems on curved domains, discretized on triangular meshes that are fitted to the boundary but unfitted to the interface. The discrete space combines GC-FE spaces on curved boundary elements, GC-IFE spaces on interface elements, and Lagrange spaces elsewhere, and discontinuous coupling is used only on the edges of boundary and interface elements, so that the number of unknowns is asymptotically that of a continuous Galerkin method. The central theoretical results are polynomial norm equivalences, inverse and trace inequalities on curved Frenet triangles, and a trace inequality for GC-IFE functions with explicit dependence on the coefficient contrast. These are employed in the first error analysis of GC-IFE spaces on triangular unfitted meshes and yield optimal-order energy-norm estimates for the proposed SDG method, with $L^2$ estimates under the standard 
dual regularity assumption. Numerical experiments confirm the predicted convergence rates with moderate and large coefficient contrasts. 

%Future work includes the proof of the low-regularity approximation property needed for 
%the $L^2$ estimate when $m\ge2$, the treatment of nonhomogeneous jump conditions 
%through enrichment, and the extension to three dimensions.

\bibliographystyle{unsrtnat}
\bibliography{FangLiuBibDesk}

\end{document}